\documentclass[12pt,reqno]{amsart}
\usepackage{amsfonts}
\usepackage{latexsym,amsmath,amsthm,amssymb,amsfonts}
\usepackage{tikz}
\numberwithin{equation}{section}
\allowdisplaybreaks

\def\CS{\mathcal S}

\def\s{\,\,\,\,}
\def\R{\mathbb{R}}

\def\Sp{\mathbb{S}}

\def\endproof{$\hfill\Box$\\}
\def\s{\,\,\,\,}
\def\R{\mathbb{R}}

\def\vol{\mbox{vol}}

\numberwithin{equation}{section}
\newtheorem{theorem}{Theorem}[section]
\newtheorem{lem}[theorem]{Lemma}
\newtheorem{thm}[theorem]{Theorem}
\newtheorem{pro}[theorem]{Proposition}
\newtheorem{cor}[theorem]{Corollary}
\newtheorem{defi}[theorem]{Definition}

\def\dint{\displaystyle{\int}}

\newcounter{Cnumber}

\def\conf{\mathcal{C}}

\title[ ]
{\bf Conformal Metric Sequences
With  Integral-bounded Scalar Curvature
}
\author[ ]
{Yuxiang Li , Zhipeng Zhou}
\address{\newline
Yuxiang Li:
 Department of Mathematical Sciences, Tsinghua University, Beijing 100084, P.R. China.
{\tt Email:yxli@math.tsinghua.edu.cn}
\newline
\newline
 Zhipeng Zhou:
Academy of Mathematics and Systems Science, CAS,
Beijing 100190, P.R. China.
{\tt Email:zhouzhipeng113@mails.ucas.ac.cn}}
\begin{document}
\maketitle

\begin{abstract}
Let $(M,g)$ be a smooth compact Riemiannian manifold
without boundary and $g_k$ be a metric conformal to $g$.  Suppose
$\vol(M,g_k)+\|R_k\|_{L^p(M,g_k)}<C$, where $R_k$ is the scalar curvature and
$p>\frac{n}{2}$.
We will use the 3-circles theorem and the John-Nirenberg inequality
to study the bubble tree convergence of $g_k$.
\end{abstract}

\section{Introduction}
The goal of this paper is to study the Gromov-Hausdorff limit of a manifold
sequence with bounded volume and $L^p$ norm of scalar curvature, in a
fixed conformal class.

Let $(M_k,g_k)$ be an $n$-dimensional Riemannian manifold. To get  a good limit, in addtion to the  curvature,
we usually assume the manifold has a lower bound on  injectivity
radius or harmonic  radius.
For example, if the section curvature $|K|\leq \Lambda$ and injectivity radius  $inj \geq i_{0}>0$, then for any $x_k\in M_k$, $(M_k,g_k,x_k)$ is pre-compact in the $C^{1,\alpha}$ topology  \cite{Cheeger}.
In \cite{Anderson}, Michael T. Anderson prove that  a similar result holds  when Ricci curvature is bounded
and the injectivity radius is bounded below.  In the same paper,  Anderson remark that a compactness result in $C^{1,\alpha}$ holds if we replace the point-wise bound on the Ricci curvature by a $L^p$-bounded with $p>\frac{n}{2}$. For  further results about the convergence under
the assumption of $\|Ric\|_{L^p}$, one can refer to \cite{Petersen-Wei,Petersen-Wei2,Adnderson-Cheeger,Cheeger2,Cheeger3,Cheeger-Colding-Tian,Tian-Zhang}.

It is natural to ask whether there is any compactness result if we only assume $L^p$-bounds on the scalar curvature instead of the Ricci curvature. From the PDE perspective, the Ricci curvature satisfies an elliptic equation (in harmonic coordinates), while there is no simple equations for the scalar curvature in general. Thus it seems necessary to impose additional assumptions other than positive injectivity radius. On the other hand, it is well-known that the scalar curvature satisfies a nice equation under conformal transformations, thus it is reasonable to replace
boundness of injectivity radius with some suitable
assumptions on the compactness of conformal classes.


More precisely, let $(M^{n},g)$ be a smooth compact Riemannian manifold of dimension $n\geq 3$ without boundary, and $R(g)$ be its scalar curvature. Denote the set of conformal metrics of $g$ as $\conf(M,g)$. For any $g'=u^\frac{4}{n-2}g\in \conf(M,g)$, the scalar curvature of $g'$ is
$$
R(g')=u^{-\frac{n+2}{n-2}}(R(g)u-4\frac{n-1}{n-2}\Delta_gu),
$$
or equivalently, $u$ satisfies the equation
\begin{equation}\label{Scalar1}
-\Delta_gu+c(n)R(g)u=c(n)R(g')u^\frac{n+2}{n-2},
\end{equation}
where $c(n)=\frac{n-2}{4(n-1)}$ and $\Delta_{g}$ is the Laplace-Beltrami operator of $g$. Sometimes, we also set $v=\log u^\frac{4}{n-2}$ and consider the following equation
\begin{equation}\label{Scalar2}
-\Delta_g v-\frac{n-2}{4}|\nabla_gv|^2+\frac{R(g)}{n-1}=\frac{R(g')}{n-1}e^v.
\end{equation}

Our goal is to study the convergence of a sequence of manifolds by using the above equations. Ideally, for a sequence of Riemannian manifolds $(M_k,g_k)$, suppose each manifold $(M_k,g_k)$ can be divided into finite many parts $(M_{kl},g_{kl})$ such that each part is almost conformal to some manifold for fixed $l$. Then we use a priori estimates of equation (\ref{Scalar1}) and (\ref{Scalar2}) to show the convergence of $(M_{kl},g_{kl})$ as $k\to\infty$. However, such convergence may have finitely many singularities. We will further use blow-up analysis and construct bubble tree convergence to study the structure at those singularities, which leads to a good limit of $(M_k,g_k)$.

As our first step towards such a program, we will fix the conformal class in this paper.
As a matter of fact, the convergence of metrics in a fixed conformal
class has been widely studied. Let us first recall some results in this field.

When $R(g_k)$ is a fixed constant, Richard Schoen raised the problem of compactness of the full set of constant scalar curvature metrics, and solved the problem for $3\leq n \leq 24$ in a joint work with M.A. Khuri and F.C. Marques~\cite{Khuri-Marques-Schoen}. The problem is also solved in some special cases, see \cite{Schoen-Zhang,Druet, Li-Zhang, Li-Zhang2, Marques}. Surprisingly, Simon Brendle constructed counterexamples of $C^{\infty}$ metrics on spheres of dimension at least 52. In a subsequent paper, Brendle and Marques extended these examples to dimension $25\leq n\leq 51$.

In \cite{Chang-Yang,Chang-Yang2} Sun-Yung A. Chang  and Paul Yang prove that when $\dim M=3$,
the space
$$
\{(M,\bar g): \bar{g}\in \conf(M,g),
\dint_{M}(1+|Ric(\bar g)|^{2})dV_{\bar g}\leq C,\lambda_{1}(\bar g)\geq \Lambda >0
\}
$$
is compact in $C^{0,\alpha}$ topology. Here $\lambda_1$ is the first
eigenvalue of $\bar{g}$.
In \cite{Gursky}, Matthew J.Gursky  proved that when $p>\frac{n}{2}$, the space
$$
\{(M,\bar g): \bar{g}\in \conf(M,g),
\dint_{M}(1+|K(\bar g)|^{p})dV_{\bar g}\leq C
\}
$$
 is compact in the $C^{0,\alpha}$ topology.\\

In this paper, we consider the convergence of a sequence $g_k=u_k^\frac{4}{n-2}g\in\conf(M,g)$ which satisfies
\begin{equation}\label{a1}
\vol(g_k)=1,\s \int_M|R(g_k)|^pdV_{g_k}<\Lambda,\s p>\frac{n}{2},
\end{equation}
where $n=\dim M>2$. It is easy to see that $u_k$ converges weakly in $W^{1,2}$ to a limit $u_0$.
But by the results of Brendle and Marques,  $g_k$ may blow up even when $R(g_k)$
is a fixed constant, thus what we are really interested in is the blowup behavior of $g_k$.

Our basic tools to study the blowup of $u_k$ are the so-called  $\epsilon$-regularity and
Three Circles Theorem.
The $\epsilon$-regularity says that if $\vol(B_r^g(p),g_k)$ is sufficiently small,
then $u_k$ is bounded in $W^{2,p}(B_{\frac{r}{2}}(p))$. Hence we can
find a finite set $\CS$ such that a subsequence of $u_k$ converges weakly
in $W^{2,p}_{loc}(M\setminus\CS)$. Then we construct a bubble tree limit near each singular point $p\in\CS$.
The Three Circles Theorem is adopted to
 investigate further details of the bubble tree convergence,
including a volume identity and a no neck result.
We remark that the usual Pohaezev inequality (see page 7 in \cite{Marques}) dose not apply here, since the scalar curvature $R$ only belongs to the $L^p$-space. Fortunately, we are able to prove a quite strong version of Three Circles Theorem, which is inspired by Qing and Tian in their study of harmonic maps \cite{Qing,Qing-Tian}. Roughly speaking, the annular area near the singularity $B_\delta(p)\setminus B_r(p)$ is almost conformal
to a long cylinder $S^{n-1}\times [-\log \delta,-\log r]$ when $\delta$ and $r$ are sufficiently small, and the equation of $u_k$
in this coordinate  is very closed to
$$
-\Delta v+\frac{(n-2)^{2}}{4}v=0,
$$
which corresponds to a strictly positive operator. Then the Three Circles Theorem asserts that the energy of $u_k$ decays exponentially as follows
$$
\int_{S^{n-1}\times [k,k+1]}u^2dtdS<C(e^{-ak}+e^{-a(-\log\delta-k)}),
$$
where $a>0$.

Our first main result is the following theorem.
\begin{thm}\label{main1}
Let $g_k=u_k^\frac{4}{n-2}g\in \mathcal{C}(M,g)$ and satisfy \eqref{a1},
and $u_0$ be the weak limit of $u_k$ in $W^{1,2}$. There
exists a finite set $\CS\subset M$, such that after passing to a
subsequence, $u_k$ converges weakly to $u_0$
in $W^{2,p}(\Omega,g)$ for any $\Omega\subset\subset M\setminus\CS$. Moreover,  $u_0$ is
a positive or equivalent to 0.
For each $x_0\in \CS$, in a local coordinate system around $x_0$,
there exists $x_k^\alpha\rightarrow x_0$, $r_k^\alpha\rightarrow 0$,  $\alpha=1$,
$\cdots$, $m$, such that

1)$$
\frac{|x_k^\alpha-x_k^\beta|}{r_k^\alpha+r_k^\beta}+\frac{r_k^\beta}{r_k^\alpha}
+\frac{r_k^\alpha}{r_k^\beta}\rightarrow +\infty,
$$

2) $(r_k^\alpha)^\frac{n-2}{2}u_k^\alpha(x_k^\alpha+r_k^\alpha x)$ converges weakly in $W^{2,p}(\R^n\setminus \CS^\alpha)$
to a positive function $ v^{\alpha}$, where $\CS^\alpha$ is a finite set.

3) $g^\alpha=v^{\alpha\frac{4}{n-2}}g_{\R^n\setminus \CS^\alpha}$ extends as
a $W^{2,p}$ metric over $S^n$.

4)
$$
 \vol(M,u_0^\frac{4}{n-2}g)+\sum_{\alpha=1}^m\vol(S^n,g^ \alpha)=1.
$$

5)  $(\cup_\alpha(S^n,g^\alpha))\cup(M,u_0^\frac{4}{n-2}g)$ is connected.\\
\end{thm}

If the weak limit $u_0=0$ and the metric degenerates away  from the singular set $\CS$, we may rescale $u_k$ and get the following result.
\begin{thm}\label{main2}
If $u_0=0$, then there exists $c_k\rightarrow+\infty$, such that
$c_ku_k$ converges weakly to a limit $G$ in $W^{1,q}(M,g)$ for any $q\in(1,\frac{n}{n-1})$, and in $W^{2,p}_{loc}(M\setminus\CS,g)$. Moreover, $G$ is positive and satisfies the equation
$$
-\Delta_g G+c(n)R(g)G=\sum_{y\in\CS}\lambda_y\delta_y,
$$
where  $\lambda_y\geq 0$ and $\delta_y$ is the Dirac function at $y$.
In particular, when $\lambda_y=0$, $G$ is smooth near
$p$; when $\lambda_y>0$  ,if we have only one bubble, then there exists a constant
$C>0$, such that
$$
\frac{1}{C}<\frac{\|u_k\|_{L^\infty(B_r(y))}}{c_k}<C,
$$
for sufficiently large $k$ and small $r$.\\
\end{thm}

In other words, Theorem~\ref{main2} says that if $u_0=0$ and
 there are exactly $m$ points $y_i$, $\cdots$, $y_m\in \CS$,  such that $\lambda_{y_i}>0$, then for any $y\in M\setminus\CS$,  $(M\setminus\CS,c_{k}^{\frac{4}{n-2}}g_{k},y)$ converges in $C^{1,\alpha}$ to the manifold $(M\setminus\{y_1,y_2,\cdots,y_m\},G^\frac{4}{n-2}g,y_\infty)$, which is  scalar flat and  complete. Also note that when $R(g)>0$, we have $\sum_{y\in\CS}\lambda_y>0$.

\begin{center}
\begin{tikzpicture}

\draw (4,-4) arc (-50: 280: 1.3);

\draw plot[smooth]coordinates{(4,-4)(3.8,-4.3)(3.6,-4.3)(3.39,-4.285)};

\draw(3.75,-4.23) ellipse (.07 and 0.03);

\draw[black!30](3.75,-4.2) ellipse (.2);

\draw plot[smooth]coordinates{(12.7,-1.6-0.8)(12,-2-.8)(11.4,-2.9-0.8)(10.8,-2.9-.8)(10.2,-2.9-.8)(9.2,-2.9-.8)};

\draw(11.1,-2.4-.69-0.4) ellipse (.28 and 0.12);

\draw[->,dashed]plot[smooth]coordinates{
(4.1,-4.1)(7,-4.05)(9,-3.8)};

\draw (7,-3.8) node {$\times c_k^\frac{4}{n-2}$};
\end{tikzpicture}

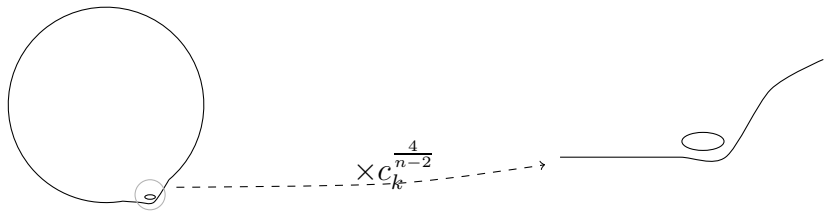
{\bf Fig 1.} When $u_0=0$ and $\lambda_y>0$, we can find $c_k$ and $y$,\\
 such that
$(M\setminus\mathcal{S},c_k^\frac{4}{n-2}g_k,y)$ converges to a\\
complete scalar flat manifold.
\end{center}

\vspace{2ex}

When $u_0>0$ or when $u_0=0$ and the limit $G$ is smooth at some singular point in $\CS$, the collar area near the singular point
shrinks to a point. Thus it is important to recover the topological information of the collar area which disappears in the limit. For this reason, we want to find a sequence of rescaling parameter $c_k\rightarrow+\infty$, such that $c_ku_k$ converges to some limit in this area.
However, if we simply multiply $u_k$ by a constant $c_k$, it is somehow difficult to control the norms of $c_ku_k$. Thus we are lead to study the function $v_k=\log u_k^{\frac{4}{n-2}}$ via equation \eqref{Scalar2} instead of \eqref{Scalar1}. We find it very convenient to work with equation \eqref{Scalar2} by applying the John-Nirenberg inequality. As a result, we get the following theorem.


\begin{thm}\label{main3}
Assume one of the following holds:

1) $u_0\neq 0$;

2) $G$ is smooth near  $y\in\CS$;

3) ${u_k}$ has at least 2 bubbles at $y\in \CS$.

Then  we can find  $y_k\rightarrow y$, $\tilde{c}_k\rightarrow+\infty$,
$t_k\rightarrow 0$, and a finite set $\tilde{\CS}\subset\R^n$,
such that $\tilde{c}_{k}t_k^\frac{n-2}{2}u_k(t_kx+y_k)$ converges weakly in $W^{2,p}_{loc}(\R^n\setminus \CS')$ to function
$v$, which is
harmonic and positive on  $\R^n\setminus\tilde{\CS}$.
Moreover, $v$ has at least 2 singular points (probably including $\infty$).\\
\end{thm}

By Bocher's Theorem (Theorem 3.9  \cite{Sheldon-Paul-Wade}), it is known that $v(x)=v_0+a|x-\tilde{y}|^{2-n}$ near $\tilde{y}\in \tilde{\CS}$,
where $v_0$ is a smooth harmonic function defined on a neighborhood of $\tilde{y}$. We say $v$ is singular at $\tilde{y}$ if $a\neq 0$, and
$v$ is singular at $\infty$, if $0$ is singular for the function $v(\frac{x}{|x|^2})$. Obviously, $v^\frac{4}{n-2}g_{\R^n}$ is smooth
across non-singular point in $\tilde{\CS}$ and  complete near each $y\in
\tilde{\CS}$. Thus, in the language of Gromov-Hausdorff limit, Theorem~\ref{main3}
can be stated as follows:
\begin{cor}\label{main3-2}
When one of 1), 2), 3) in Theorem \ref{main3} holds, we can find
a  closed neighborhood $\Omega_k$ of $y$,   and $y_k\in\Omega_k$, $\hat{c}_k\rightarrow+\infty$, such that $(\Omega_k,\hat{c}_kg_k,y_k)$ converges
to $(M_\infty,g_\infty,y_\infty)$ in Gromov-Hausdorff limit, where $(M_\infty,g_\infty,y_\infty)$
is complete, smooth, scalar flat, and conformally diffeomeomorphic to a $n$-sphere
with at least 2 points removed.\\
\end{cor}

\begin{center}
\begin{tikzpicture}

\draw (4,-4) arc (-66: 280: 1.3);
\draw (3.8,-4.25) arc (120: 450: .43 and .33);
\draw plot[smooth,tension=1.9]coordinates{(4,-4)(3.9,-4.1)(4.05,-4.2)};
\draw plot[smooth,tension=1.9]coordinates{(3.8,-4.25)(3.85,-4.15)(3.7,-4.1)};

\draw(4.2,-4.5) ellipse (.15 and 0.1);

\draw[blue!30](3.9,-4.15) ellipse (.15);

\draw (10,-4.5) arc (-50: 50: 2);
\draw (12,-1.4) arc (130: 230: 2);

\draw[->,dashed]plot[smooth]coordinates{
(4.1,-4.1)(7,-3.75)(10,-3)};

\draw (7,-3.4) node {$\times \hat{c}_k$};
\end{tikzpicture}

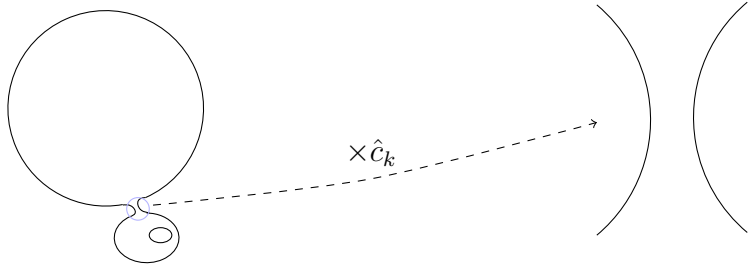
{\bf Fig 2.} When one of 1), 2), 3) in Theorem \ref{main3} holds, we can\\ find
$\Omega_k$ of $y$, $y_k$, $\hat{c}_k\rightarrow+\infty$, such that
$(\Omega_k,\hat{c}_kg_k,y_k)$\\ converges to a complete scalar flat manifold.
\end{center}
\vspace{2ex}

Now we discuss some simple applications of our results. First note that a simple connected flat manifold with an end collared topologically by $S^{n-1}\times\R$
must be $\R^n$. Thus if we replace the assumption $\|R\|_{L^p}<C$ by $L^p$-boundedness of the sectional curvature $\|K\|_{L^p}<C$, the above theorems will imply
Gursky's result \cite{Gursky}.

Another quick application is to prove the compactness of a sequence of manifolds with $L^p$-bounded Ricci curvature.
\begin{cor}\label{main4}
Let $g_k=u_k^\frac{4}{n-2}g$, which
$$
\vol(M,g_k)=1,\s \int_M|Ric(g_k)|^pdV_{g_k}<\Lambda,\s p>\frac{n}{2}.
$$
Then the Gromov-Hausdorff limit of $(M,g_k)$ is either $M$ or $S^{n}$, endowed with
a $W^{2,p}$-metric.\\
\end{cor}

The rest of the paper is organized as follows.
In Section 2, we list some preliminary results in elliptic PDE theory which are needed in the paper.
In Section 3, we derive a gap theorem and a non-collapsing theorem.
In Section 4 ,we prove the powerful Three Circles Theorem.
In Section 5, we construct the bubble tree and prove Theorem~\ref{main1}. Theorem \ref{main2} and \ref{main3} are proved in Section 6 and Section 7 respectively.\\

{\bf Acknowledgement:} The authors would like to thank Prof. Hao Yin for bringing Three Circles Theorem to our attention. The authors also thank Prof. Chong Song for helpful suggestions during the preparation of this paper.

\section{Preliminary}
In this section, we list some results on elliptic equation
which will be used later.
\begin{lem}\label{epsilon}
 Let $B$ be the unit ball in $\R^n$. $u\in W^{2,p}(B)$ is a solution to the following equation with $p > \frac{n}{2}$:
\begin{equation}\label{equation.epsilon}
a^{ij}u_{ij}+b_iu_i+cu=R |u|^\frac{4}{n-2}u+f(x)
\end{equation}
where $f\in L^p(B)$ and $R$ is a measurable function with
$$
\int_B|R|^p |u|^\frac{2n}{n-2}dx<\Lambda.
$$
We assume
$$
\lambda_1<a^{ij}<\Lambda_1,\s \|a^{ij}\|_{C^{0,1}(B)}+\|b^i\|_{C^{0,1}(B)}+
\|c\|_{C^{0,1}(B)}<\Lambda_2.
$$
Then there exists an $\epsilon_0=\epsilon_0(\Lambda,\Lambda_1,\lambda_1,\Lambda_2)>0$, such that if $\int_B |u|^\frac{2n}{n-2}dx<\epsilon_0$,
then
$$
\|u\|_{W^{2,p}(B_\frac{1}{2})}<C(\|u\|_{L^\frac{2n}{n-2}(B)}+\|f\|_{L^p(B)}),
$$
where $C$ only depends on $\Lambda$, $\Lambda_1$, $\lambda_1$, $\Lambda_2$ and $\epsilon_0$.
\end{lem}

\proof

Let $\eta$ be a smooth function. We have
\begin{eqnarray}\label{Lgamma}
\left\|R\eta u^{\frac{n+2}{n-2}}\right\|_{L^{\gamma}(B)}
&=&\left(\dint_{B(p)}|R|^{\gamma}\eta^{\gamma}  |u|^{\frac{n+2}{n-2}\gamma}dx\right)^{\frac{1}{\gamma}}\\\nonumber
&=&\left(\dint_{B}|R|^{\gamma}{|u|}^{\frac{2n}{n-2}\frac{2\gamma}{n}}{|u|}^{\frac{n+2}{n-2}\gamma-\frac{2n}{n-2}\frac{2\gamma}{n}}\eta^{\gamma}dx\right)^{\frac{1}{\gamma}}\\\nonumber
&\leq&\left(\dint_{B}|R|^{\frac{n}{2}}{ |u|}^{\frac{2n}{n-2}}dx\right)^{\frac{2}{n}}\left(\dint_{B}({ |u|}\eta)^{\frac{n\gamma}{n-2\gamma}}dx\right)^{\frac{1}{\gamma}-\frac{2}{n}}\\\nonumber
&\leq&\left(\dint_{B}|R|^{p}{ |u|}^{\frac{2n}{n-2}}dx\right)^{\frac{1}{p}}\left(\dint_{B}{|u|}^{\frac{2n}{n-2}}dx\right)^{\frac{2}{n}-\frac{1}{p}}\left(\dint_{B}(\eta { |u|})^{\frac{n\gamma}{n-2\gamma}}dx\right)^{\frac{n-2\gamma}{n\gamma}}
\end{eqnarray}

Note that when
$\gamma=\frac{2n}{n+2}$, $\frac{n\gamma}{n-2\gamma}=\frac{2n}{n-2}$. Putting $\eta=1$, we have
$$
\|Ru^\frac{n+2}{n-2}\|_{L^\frac{2n}{n+2}(B)}<C\|u\|_{L^\frac{2n}{n-2}(B)}.
$$
Here C depends on the $\epsilon_0$ and $\Lambda$.
Thus we get
$$
\|u\|_{W^{2,\frac{2n}{n+2}}(B_r)}<C(r,\Lambda,\epsilon_{0})(\|f\|_{L^\frac{2n}{n+2}(B)}+\|u\|_{L^\frac{2n}{n-2}(B)}),
$$
where $r\in(0,1)$.
By Sobolev embedding, we get
$$
\|\nabla u\|_{L^2(B_r)}<C(\|f\|_{L^\frac{2n}{n+2}(B)}+\|u\|_{L^\frac{2n}{n-2}(B)})<C(\|f\|_{L^p(B)}+
\|u\|_{L^\frac{2n}{n-2}(B)}).
$$

Now,
let $\eta$ be a cut-off function which is $1$ in $B_{\frac{1}{2}}$ and $0$ in $ B_{\frac{2}{3}}^c$. We have
\begin{eqnarray*}
a^{ij} (\eta u)_{ij}&=&
a^{ij}\eta_{ij} u+2a^{ij}u_i\eta_j+\eta a^{ij}u_{ij}\\
&=&a^{ij}\eta_{ij} u+2a^{ij}u_i\eta_j
-\eta b_i u_i-c\eta u+\eta R|u|^\frac{4}{n-2}u+\eta f(x)\\
&:=&\eta R{ |u|}^\frac{4}{n-2}u+F,
\end{eqnarray*}
where
$$
|F|\leq C(|\nabla u|+|u|+|f|).
$$

By \eqref{Lgamma},
$$
\left\|\eta Ru^\frac{n+2}{n-2}\right\|_{L^\gamma(B)}
\leq C\epsilon_0^{\frac{2}{n}-\frac{1}{p}}
\|\eta u\|_{W^{2,\gamma}(B)}.
$$
Let $\frac{2n}{n+2}<\gamma_0<2$, and $\epsilon_0$ be sufficiently small. By elliptic estimate and H\"older
inequality, we get
$$
\|\eta u\|_{W^{2,\gamma_0}(B)}\leq
C\|F\|_{L^{\gamma_0}(B)}<C(\|u\|_{L^\frac{2n}{n-2}(B)}+\|f\|_{L^p(B)}).
$$

Then
\begin{eqnarray*}
\left\|R u^{\frac{n+2}{n-2}}\right\|_{L^{\gamma}(B_r)}
&=&\left(\dint_{B_r}|R|^{\gamma}{ |u|}^{\frac{2n}{n-2}\frac{\gamma}{p}}{ |u|}^{\frac{n+2}{n-2}\gamma-\frac{2n}{n-2}\frac{\gamma}{p}}dx\right)^{\frac{1}{\gamma}}\\
&\leq&\left(\dint_{B_r}|R|^{p}{  |u|}^{\frac{2n}{n-2}}dx\right)^{\frac{1}{p}}\left(\dint_{B_r}{  |u|}^{(\frac{n+2}{n-2}\gamma-\frac{2n}{n-2}\frac{\gamma}{p})\frac{p}{p-\gamma}}dx\right)^{\frac{p-\gamma}{p\gamma}}.
\end{eqnarray*}
We set
$$
(\frac{n+2}{n-2}\gamma_{k}-\frac{2n}{n-2}\frac{\gamma_{k}}{p})
\frac{p}{p-\gamma_{k}}=s_{k},\s
s_{k}=\frac{n\gamma_{k-1}}{n-2\gamma_{k-1}},\s k\geq 1.
$$
Thus, if $s_k$ and $\gamma_k$ make sense, we have
$$
\|R u^{\frac{n+2}{n-2}}\|_{L^{\gamma_k}(B_r)}\leq
C\|u\|_{L^{s_k}(B)}^\frac{(p-\gamma_k)s_k}{p\gamma_k},
$$
which implies an estimate of $\|u\|_{L^{s_{k+1}}}$.

We claim that, there exists $k$   so that  $\gamma_k\geq
\frac{n}{2}$. Assume this is not true.
Let $\tau=\frac{n+2}{n-2}-\frac{2n}{n-2}\frac{1}{p}$ which
is greater than 1. We have
$$
\frac{1}{\tau}(\frac{1}{s_{k+1}}+\frac{2}{n}-\frac{1}{p})=\frac{1}{s_{k}},
$$
then
$$
\frac{1}{\tau^{k+1}}\frac{1}{s_{k+1}}-\frac{1}{\tau^{k}}\frac{1}{s_{k}}=-\frac{1}{\tau^{k+1}}(\frac{2}{n}-\frac{1}{p}).
$$
We get
\begin{eqnarray*}
\lim_{k\to\infty}{\frac{1}{\tau^{k+1}s_{k+1}}}&=&\frac{1}{s_{0}}-\sum_{i=1}^{\infty}\frac{1}{\tau^{i}}(\frac{2}{n}-\frac{1}{p})\\
&=&\frac{1}{s_0}-\frac{1}{\tau-1}(\frac{2}{n}-\frac{1}{p})\\
&=&\frac{1}{s_0}-\frac{n-2}{2n}\\
&<&0,
\end{eqnarray*}
which contradicts   the  assumption that   $\gamma_{k}<\frac{n}{2}$.

Now, we can assume  $\gamma_k\geq\frac{n}{2}$. By Sobolev
embedding, $u$ is bounded in $L^q(B_r)$ for any $q\geq 0$
and  $r\in(0,1)$. Then $Ru^\frac{n+2}{n-2}$  is bounded in
 $L^p(B_{r})$. We complete the  proof of this theorem.
~\endproof

\begin{lem}\label{A1}
Let $\Omega$ be a bounded domain in $\R^n$, and
$v$ solve the equation:
$$
-\Delta v+a(x)v=R|u|^\frac{4}{n-2}v, \s x\in\Omega.
$$
We assume $\|a\|_{C^0}<A_0$, $\int_\Omega|R|^p|u|^\frac{2n}{n-2}<A_1$, and
$\|u\|_{L^\infty(\Omega)}<A_2$, where $p>\frac{n}{2}$. Then for any $q>\frac{p}{p-1}$, and $\Omega'\subset\subset\Omega$,
we have
$$
\|v\|_{W^{2,p}(\Omega')}\leq C(A_0,A_1,A_2,p,q,\Omega,\Omega')\|v\|_{L^q(\Omega)}.
$$
\end{lem}

\proof

For any $s<p$, we have
\begin{eqnarray*}
\int|R|^s|u|^{\frac{4}{n-2}s}|v|^s
&=&\int|R|^s|u|^{\frac{2n}{n-2}\frac{s}{p}}
|u|^{\frac{4s}{n-2}(1-\frac{n}{2p})}|v|^s\\
&\leq&\|u\|_{L^\infty}^{\frac{4s}{n-2}(1-\frac{n}{2p})}\left(\int|R|^p|u|^\frac{2n}{n-2}\right)^\frac{s}{p}
\left(\int |v|^\frac{ps}{p-s}\right)^\frac{p-s}{p}.
\end{eqnarray*}
Let  $\frac{ps_1}{p-s_1}=q$, i.e.
$$
\frac{1}{s_1}=\frac{1}{p}+\frac{1}{q}.
$$
Since $q>\frac{p}{p-1}$, $s_1>1$. Then we get the estimate of
$W^{2,s_1}$-norm of $v$, and hence the $L^\frac{ns_1}{n-2s_1}$-norm of
$v$.
Let $q_1=q$ and $q_2=\frac{ns_1}{n-2s_1}$.
We have
$$
\frac{1}{q_2}=\frac{1}{q_1}+\frac{1}{p}-\frac{2}{n}.
$$
In the same way, we  get the estimate of
$\|v\|_{L^{s_i}}$, where $s_i$ is the solution to
$q_i=\frac{ps_i}{p-s_i}$ and
$
\frac{1}{q_i}=\frac{1}{q_{i-1}}+\frac{1}{p}-\frac{2}{n}$.

Note that $\frac{1}{p}-\frac{2}{n}<0$, we can get $i$,
such that $s_i\geq \frac{n}{2}$. Then we get estimate of
$\|v\|_{L^\infty}$, and then the estimate of
$\|v\|_{W^{2,p}}$.
\endproof

\begin{lem}\label{Lq}
Let $u\in W^{1,2}(M,g)$ solve the equation
$-\Delta u+u=f$, where  $\|f\|_{L^1(M)=1}$.
Then for any $q_1\in(1,\frac{n}{n-2})$, $q_2\in(1,\frac{n}{n-1})$,
we have
$$
\int_{B_r(x)}|u|^{q_1}\leq C(M,g,q_1)r^{n-(n-2)q_1},\s
\int_{B_r(x)}|\nabla u|^{q_2}<C(M,g,q_2)r^{n-(n-1)q_2}.
$$
\end{lem}
\proof
Let $G(x,y)$ be the Green function of $-\Delta_{g}+1$. Then we have $G=r^{2-n}(1+\psi_{1})$, $|\nabla G|\leq cr^{1-n}(1+\psi_{2})$, where
$\psi_{1}=o(1)$ and $\psi_{2}=o(1)$. Moreover, we have
$$
u(x)=\dint_{M}G(x,y)f(y)dV_{g}(y).
$$
Then for any $q_1\in(1,\frac{n}{n-2})$, by Jensen inequality, we  get
\begin{eqnarray*}
\int_{B_r(z)}|u|^{q_1}(y)dV_g(y)&=&
\int_{B_{r}(z)}\left|\int_{M}G(x,y)f(y)dV_{g}(y)\right|^{q_1}dV_{g}(x)\\
&\leq&\int_{B_{r}(z)}\left(\int_{M}|G(x,y)||f(y)|dV_{g}(y)\right)^{q_1}dV_{g}(x)\\
&\leq&\dint_{B_{r}(z)}\dint_{M}|G(x,y)|^{q_1}|f(y)|dV_{g}(y) dV_{g}(x)\\
&\leq& \dint_{B_{r}(z)}(r^{2-n}(1+\psi_{1}))^{q_1}dV_{g}(x)\\
&\leq&{ C(M,g,q_{1})}r^{n-(n-2)q_1}.
\end{eqnarray*}

Next, we estimate $\|\nabla u\|_{L^{q_2}}$. Since
$$
\nabla u(x)=\dint_{M}\nabla G(x,y)f(y)dV_{g}(y),
$$
for any $q_2\in(1,\frac{n}{n-1})$,
we have
\begin{eqnarray*}
\int_{B_{r}(z)}|\nabla u(x)|^{q_2}dV_{g}(x)&\leq
&\int_{B_{r}(z)}\int_{M}|\nabla G(x,y)|f(y)dV_{g}(y))^{q_2}dV_{g}(x)\\
&=&\int_{B_{r}(z)}(\int_{M}\nabla G(x,y)|f(y)|dV_{g}(y))^{q_2}dV_{g}(x)\\
&\leq&\int_{B_{r}(z)}(\dint_{M}|\nabla G(x,y)|^{q_2}|f(y)|dV_{g}(y))dV_{g}(x)\\
&\leq& \dint_{B_{r}(z)}(r^{1-n}(1+\psi_{2}))^{q}dV_{g}(x)\\
&\leq&C(M,g,q_2)r^{n-(n-1)q_2}.
\end{eqnarray*}
\endproof

For a proof of the following lemma, one can refers to  \cite{Lin-Han}

\begin{lem}\label{moser}
Let $u$ be a non-negative function in $W^{1,2}(B)$ which
satisfies
$$
-D_{j}(a^{ij}D_{i}u)\leq f.
$$
We assume
$$
0<\lambda< a^{ij}(x)<\Lambda,
$$
Then
$$
\|u\|_{L^\infty(B_\frac{1}{2})}\leq C(p,q)(\|u\|_{L^q}+\|f\|_{L^p}),
$$
where $p>\frac{n}{2}$ and $q>0$.
\end{lem}

We also need the following  version of John-Nirenberg inequality:
\begin{thm}\label{J-N}
Let $B_1$ be the unit ball and $u\in W^{1,2}(B_1)$. There exist
positive constants $\alpha$ and $\beta$, such that if
$$
\int_{B_r(x)}|\nabla u|^2dx<r^{n-2},\s \forall B_r(x)\subset B_1,
$$
then
$$
\int_{B_1}e^{\alpha u}dx\int_{B_1}e^{-\alpha u}dx<\beta.
$$
\end{thm}

As an application, we have the following:

\begin{cor}\label{positive}
Let $u$ be a non-negative function in $W^{1,2}\cap C^0(B)$ which
satisfies
$$
-D_{j}(a^{ij}D_{i}u)=au+fu^\frac{n+2}{n-2}.
$$
We assume
$$
0<\lambda< a^{ij}(x)<\Lambda,\s a\in C^0(B)
$$
and
$$
\int_B|f|^{p}u^\frac{2n}{n-2}<+\infty
$$
for some $p>\frac{n}{2}$.
Then either $u=0$ in $B_\frac{1}{2}$, or  $essinf_{B_\frac{1}{4}} u>0$.
\end{cor}
\proof We assume $u\neq 0$.

Let $\eta$ be a cut-off function, which is 1 on $B_{\frac{1}{2}}(0)$ and 0 outside $B_{1}(0)$. Let $\eta_r=\eta(\frac{x}{r})$.
We have
$$
\int_{B_{r}(x)}a^{ij} \nabla_{i}u\nabla_{j}(\eta_{r}^{2}\frac{1}{u+\epsilon})dx=
\int_{B_r(x)}\left(a\frac{u}{u+\epsilon}\eta_{r}^{2}+f\frac{u^{\frac{n+2}{n-2}}}{u+\epsilon}\eta_{r}^{2}\right)dx.
$$
Then
$$
\int_{B_{r}(x)}\eta_{r}^{2}|\nabla\log(u+\epsilon)|^2dx\leq C\left(\int_{B_{r}(x)}|\nabla\eta_{r}|^{2}dx+\int_{B_r(x)}(1+|f|u^{\frac{4}{n-2}})dx\right).
$$
By Sobolev embedding,  $\int_Bu^\frac{2n}{n-2}<+\infty$.
Noting that
$$
\int_{B_{r}}|f|u^\frac{4}{n-2}<C(\int_{B_{r}}|f|^pu^\frac{2n}{n-2})^\frac{1}{p}
(\int_{B_{r}}u^\frac{2n}{n-2})^{\frac{2}{n}-\frac{1}{p}}(\int_{B_{r}})^{1-\frac{2}{n}},
$$
we get
$$
\int_{B_{\frac{r}{2}}(x)}|\nabla\log(u+\epsilon)|^2\leq C(\Lambda,\lambda)r^{n-2}.
$$

Then, by John-Nirenberg inequality,
$$
\int_{B_{\frac{1}{2}}(0)}(u+\epsilon)^\frac{\alpha}{\sqrt{C(\Lambda,\lambda)}}\int_{B_{\frac{1}{2}}(0)}(u+
\epsilon)^{-\frac{\alpha}{\sqrt{C(\Lambda,\lambda)}}}<\beta.
$$
For convinience, let $\alpha'=\frac{\alpha}{\sqrt{C(\Lambda,\lambda)}}$.
Since $u\neq 0$, we can find $\delta_0$, such that
$\int_{\{u>\delta_0\}\cap B_{\frac{1}{2}}}>b>0$, which yields $\int_{B_{\frac{1}{2}}(0)}(u+\epsilon)^{\alpha'}>b\delta_0^{\alpha'}$. Thus
$$
\int_{B_{\frac{1}{2}}(0)}(u+
\epsilon)^{-\alpha'}<C_{1}.
$$
Choose  $\tau$ be sufficiently small, such that
$$
\frac{2n}{p(n-2)}<\frac{4}{n-2}-\tau,\s \tau p<\alpha'.
$$
Then
\begin{eqnarray*}
-\nabla_{j}(a^{ij}\nabla_{i}(u+\epsilon)^{-\tau})
&\leq& |a|(u+\epsilon)^{-\tau} +|f|u^{\frac{4}{n-2}-\tau}\\
&\leq& |a|(u+\epsilon)^{-\tau}+
|f|u^\frac{2n}{p(n-2)}u^{\frac{4}{n-2}-\tau-\frac{2n}{p(n-2)}}
\end{eqnarray*}
in the weak sense.
Since
$$
\int_{B_{\frac{1}{2}}(0)}|fu^\frac{2n}{p(n-2)}u^{\frac{4}{n-2}-\tau-\frac{2n}{p(n-2)}}|^p
\leq C\|u\|_{L^\infty}^{(\frac{4}{n-2}-\tau-\frac{2n}{p(n-2)})p}
\int_{B_{\frac{1}{{ 2}}}(0)}|f|^pu^\frac{2n}{n-2},
$$
by Lemma \ref{moser},  $\|(u+\epsilon)^{-\tau}\|_{L^\infty{ (B_\frac{1}{4})}}
<C$. Letting $\epsilon\rightarrow 0$, we get $\|u^{-\tau}\|_{L^\infty({  B_\frac{1}{4})}}
<C$.
\endproof

\section{ $W^{2,p}$-conformal metrics}
Let $p>\frac{n}{2}$ and $u$ be a positive function in $W^{2,p}(M,g)$.
$g'=u^\frac{4}{n-2}g$ is a $W^{2,p}$-metric,  curvature of $g'$
is well-defined.
For convenience, we denote
$$
SC^p(M,g)=\{g'=u^\frac{4}{n-2}g:u\in
W^{2,p}(M,g),\s u>0\},
$$
$$
SC^p_\Lambda(M,g)=\{g'=u^\frac{4}{n-2}g\in SC^p(M,g):\|R(g')\|^{p}_{L^p(M,g')}
\leq\Lambda\}.
$$

We have the following gap lemma:

\begin{thm}\label{gap}
Let $g'\in SC^p_\Lambda(\mathbb{S}^n,g_{\mathbb{S}^n})$. Then there exists an
$\tau=\tau(\Lambda)$, such that
$$
\vol(\mathbb{S}^n,g')\geq\tau.
$$
\end{thm}

\proof

Let $g'=u^\frac{4}{n-2}g_{\mathbb{S}^n}$.
By \eqref{Scalar1},
\begin{eqnarray*}
\int_{\mathbb{S}^n}(|\nabla_{g_{\mathbb{S}^n}}u|^2+c(n)R_{\mathbb{S}^n}u^2)dV_{g_{\mathbb{S}^n}}&=&
\int_{\mathbb{S}^n}c(n)R(g')u^\frac{2n}{n-2}dV_{\mathbb{S}^n}\\
&=&
c(n)\int_{\mathbb{S}^n}R(g')dV_{g'}\\
&\leq&
c(n)\Lambda^\frac{1}{p}(\vol(\mathbb{S}^n,g'))^\frac{p-1}{p}.
\end{eqnarray*}
Then
$$
\lambda \vol^\frac{n-2}{n}(\mathbb{S}^n,g')\leq c(n)\Lambda^\frac{1}{p}
vol^\frac{p-1}{p}(\mathbb{S}^n,g'),
$$
where $\lambda$ is the Yamabe constant of $\mathbb{S}^n$.
Since $\frac{p-1}{p}>\frac{n-2}{n}$,   we finish the proof.
\endproof

Combining Proposition \ref{epsilon} and \ref{positive} together, we get
the following:

\begin{pro}Let $B$ be the unit ball of $\R^n$ and
$\hat{g}_k$ be a smooth metric over $B$ which converges to
a metric $g$ smoothly. Let
 $g_k=u_k^\frac{4}{n-2}\hat{g}_k\in SC^p_{\Lambda}(B,\hat{g}_k)$  for some $p>\frac{n}{2}$. Then there exists an $\epsilon_0$, such that
if $\vol(B,g_k)<
\epsilon_0$, then
$$
\|u_k\|_{W^{2,p}(B_\frac{1}{2},g)}<C.
$$
Moreover, a subsequence of $\{u_k\}$ converges weakly in $W^{2,p}(B_\frac{1}{2})$. The limit is positive or equivalent to
0.
\end{pro}

Note  that in  the above proposition, the Gromov-Hausdorff limit of $(B,g_k)$  may be a point. The following result gives a sufficient condition of  non-collasping.

\begin{pro}
Let $B$, $\hat{g}$, $g_k$, $\epsilon_0$ be as in the above   proposition. If
$vol(B,g_k)<\epsilon_0$, and
$$
\inf_kd_{g_k}(0,\partial B)>0,
$$
then a subsequence of $\{g_k\}$ converges weakly in $W^{2,p}(B_\frac{1}{2})$
to a metric $g'\in SC^p(g,B_\frac{1}{2})$.
\end{pro}

\proof Without loss of generality, we assume $g_{ij}=\delta_{ij}$.
Assume $u_k\rightarrow 0$ uniformly in $B_a$
for any $a<1$, and
$$
\inf_{k}d_{g_k}(0,\partial B)>\delta>0.
$$

For any $y\in\partial B$,
$$
d_{g_k}(0,y)\leq C\int_0^{1}u^\frac{2}{n-2}(ty)dt.
$$
Then we have
\begin{eqnarray*}
c\delta&\leq&\int_{S^{n-1}}\int_0^{1}u_k^\frac{2}{n-2}drdS\\
&\leq&
\int_{S^{n-1}}\int_0^au_k^\frac{2}{n-2}drdS+
\int_{S^{n-1}}\int_a^{1}u_k^\frac{2}{n-2}r^\frac{n-1}{2}
r^\frac{-(n-1)}{2}drdS\\
&\leq&\int_{S^{n-1}}\int_0^au_k^\frac{2}{n-2}drdS+
\sqrt{\int_{S^{n-1}}\int_a^{1}u_k^\frac{4}{n-2}r^{n-1}drdS}
\sqrt{\int_{S^{n-1}}\int_a^{1}r^{1-n}drdS}\\
&\leq&\int_{S^{n-1}}\int_0^au_k^\frac{2}{n-2}drdS
+{  C(n)}(a^{2-n}-1)^{\frac{1}{2}}.
\end{eqnarray*}
We get a contradiction if we choose $a$ to be sufficiently
close to $1$.
\endproof

\section{Three Circles Theorem}
In this section, we will prove the Three Circles Theorem and give
some applications. Since it is convenient to state and prove this theorem
on pipes, we let $Q=[0,3L]\times S^{n-1}$, and
$$
Q_i=[(i-1)L,iL]\times S^{n-1},\s i=1,2,3.
$$
Set $g_0=dt^2+g_{\mathbb{S}^{n-1}}$ and $dV_0=
dV_{g_0}$. We have the following:
\begin{thm}\label{3-circle}
Let $g$ be a metric over $Q$ and $u\in W^{2,p}$
which solves the equation
$$
-\Delta_gu+c(n)R(g)u=c(n)R|u|^\frac{4}{n-2}u,
$$
where $\int_Q|R|^p|u|^\frac{2n}{n-2}dV_g<\Lambda$
 for some $p>\frac{n}{2}$.
There exist $\epsilon_1$, $\tau$, and $L_0$, such that
if
$$
\|g-g_0\|_{C^2(Q)}<\tau,\s \int_Q|u|^\frac{2n}{n-2}dV_g<\epsilon_1,
$$
then for any $L>L_0$, we have

\begin{itemize}
\item[1) ]$\dint_{Q_1}u^2dV_{0}\leq e^{-L}
\dint_{Q_2}u^2dV_{0}$ implies
$$\int_{Q_2}u^2dV_{0}\leq e^{-L}\int_{Q_3}u^2dV_{0}.$$

\item[2) ] $\dint_{Q_3}u^2dV_{0}\leq e^{-L}
\dint_{Q_2}u^2dV_{0}$ implies
$$\int_{Q_2}u^2dV_{0}\leq e^{-L}\int_{Q_1}u^2dV_{0}.$$

\item[3) ]either
$$
\int_{Q_2}u^2dV_{0} \leq e^{-L}\int_{Q_1}u^2dV_{0}\s
or\s
\int_{Q_2}u^2dV_{0}\leq e^{-L}\int_{Q_3}u^2dV_{0}.
$$
\end{itemize}
\end{thm}

We first prove this theorem for case of $g=g_0$ and $R=(n-1)(n-2)$:
\begin{lem}\label{3-circle.R=0}
Let $u\neq 0$ solve the following equation on $Q$:
$$
-\Delta u+\frac{(n-2)^{2}}{4}u=0,
$$
then there exists $L_0$, such that for any $L>L_0$, we have

\begin{itemize}
\item[1) ]$\dint_{Q_1}u^2dV_{0}\leq e^{-L}
\dint_{Q_2}u^2dV_{0}$ implies
$$\int_{Q_2}u^2dV_{0}< e^{-L}\int_{Q_3}u^2dV_{0}.$$

\item[2) ] $\dint_{Q_3}u^2dV_{0}\leq e^{-L}
\dint_{Q_2}u^2dV_{0}$ implies
$$\int_{Q_2}u^2dV_{0}< e^{-L}\int_{Q_1}u^2dV_{0}.$$

\item[3) ]either
$$
\int_{Q_2}u^2dV_{0} < e^{-L}\int_{Q_1}u^2dV_{0}\s
or\s
\int_{Q_2}u^2dV_{0}< e^{-L}\int_{Q_3}u^2dV_{0}.
$$
\end{itemize}
\end{lem}

\proof

Let  $\{\varphi_m:m=0,1,2,\cdots\}$,
be an orthonormal eigenbasis   for  $\Delta_{g_{\mathbb{S}^{n-1}}}$.
We set
$$
-\Delta \varphi_m=\lambda_m\varphi_m.
$$
It is well-known that (cf \cite{Sheldon-Paul-Wade})
$$
\lambda_m\in\{j(j+n-2):j=0,1,\cdots\}.
$$
If we assume $u=\sum\limits_{m=0}^\infty v_m(t)\varphi_m(x)$, then we have
$$
0=-v_m''(t)+(\lambda_m+\frac{(n-2)^2}{4})v_m.
$$
We get
$$
v_m=a_me^{-\tau_mt}+b_me^{\tau_mt},\s \forall m\geq 0,
$$
where $a_m$ and $b_m$ are constants, and
$$
\tau_m=\sqrt{\lambda_m+\frac{(n-2)^2}{4}}.
$$
Then we have
$$
u=\sum_{m=0}^\infty (a_m e^{-\tau_mt}+b_m e^{\tau_m t})\varphi_m.
$$
We have
\begin{eqnarray*}
\|u\|^{2}_{L^{2}([(i-1)L,iL]\times S^{n-1})}
&=&\dint_{[(i-1)L,iL]\times S^{n-1}}(\sum\limits_{m=0}^{+\infty}(a_m e^{-\tau_mt}+b_m e^{\tau_m t})\varphi_{m}(\theta))^{2}dV_0\\
&=&\dint_{(i-1)L}^{iL}\sum\limits_{m=0}^{+\infty}(a_{m}e^{-\tau_m t}+b_me^{\tau_{m}t})^2dt\\
&=&\sum\limits_{m=0}^\infty\left(
a_m^2\frac{1-e^{-2\tau_m L}}{2\tau_m}e^{-2\tau_m(i-1)L}
+b_m^2\frac{e^{2\tau_m L}-1}{2\tau_m}e^{2\tau_m(i-1)L}+2a_m b_mL\right)
\end{eqnarray*}
When $a_mb_m\leq 0$, we have
\begin{eqnarray*}
a_m^2\frac{1-e^{-2\tau_m L}}{2\tau_m}e^{-2\tau_m(2-1)L}
+b_m^2\frac{e^{2\tau_m L}-1}{2\tau_m}e^{2\tau_m(2-1)L}+2a_m b_mL\\
\leq
e^{-2\tau_mL}(a_m^2\frac{1-e^{-2\tau_m L}}{2\tau_m}+2a_m b_mL)
+e^{-2\tau_mL}(b_m^2\frac{e^{2\tau_m L}-1}{2\tau_m}e^{2\tau_m(3-1)L}
+2a_m b_mL),
\end{eqnarray*}
when $L$ is sufficiently large.
When $a_mb_m>0$, we can choose $L$ to be sufficiently large, such that
\begin{eqnarray*}
a_m^2\frac{1-e^{-2\tau_m L}}{2\tau_m}e^{-2\tau_m(2-1)L}
+b_m^2\frac{e^{2\tau_m L}-1}{2\tau_m}e^{2\tau_m(2-1)L}+2a_m b_mL\\
\leq 2a_m^2\frac{1-e^{-2\tau_m L}}{2\tau_m}e^{-2\tau_m(2-1)L}
+b_m^2\frac{2\tau_m}{1-e^{-2\tau_m L}}e^{2\tau_m(2-1)L}L^{2}
+b_m^2\frac{e^{2\tau_m L}-1}{2\tau_m}e^{2\tau_m(2-1)L}\\
\leq
3e^{-2\tau_mL}(a_m^2\frac{1-e^{-2\tau_m L}}{2\tau_m}+2a_m b_mL)
+3e^{-2\tau_mL}(b_m^2\frac{e^{2\tau_m L}-1}{2\tau_m}e^{2\tau_m(3-1)L}
+2a_m b_mL).
\end{eqnarray*}
Hence,
$$
\int_{[L,2L]\times S^{n-1}}|u|^2dV_0\leq 3e^{-\tau_0 L}(
\int_{[0,L]\times S^{n-1}}|u|^2dV_0+\int_{[2L,3L]\times S^{n-1}}|u|^2dV_0),
$$
which yields the lemma if  $L$ is chosen to be sufficiently large.
\endproof

{\it Proof of Theorem \ref{3-circle}:}
If the statement in 1) was false, then we can find $g_k$  and
$u_k$,  s.t.
$$
g_k\rightarrow g_0 \mbox{ in }C^2,\s \int_{Q}{ |u_k|}^\frac{2n}{n-2}
\rightarrow 0,$$
$$
-\Delta_{g_k}u_k+c(n)R(g_k)u_k=c(n)R_k{ |u_k|}^\frac{4}{n-2}u_{k},
$$
$$
\int_Q|R_k|^p{ |u_k|}^\frac{2n}{n-2}dV_{g_k}\leq \Lambda
$$
and
$$
\int_{Q_1}|u_k|^2dV_0\leq e^{-L}\int_{Q_2}|u_k|^2dV_0,\s
\int_{Q_2}|u_k|^2dV_0 > e^{-L}\int_{Q_3}
|u_k|^2dV_0.
$$
Let $$v_k=\frac{u_k}{\|u_k\|_{L^2(Q_2)}}.$$
We have
$$
\int_{Q_1}|v_k|^2dV_0\leq e^{-L}\int_{Q_2}|v_k|^2dV_0,\s
\int_{Q_3}|v_k|^2dV_0 < e^{L}\int_{Q_2}|v_k|^2dV_0,
$$
and
$$
\int_{Q_2}|v_k|^2dV_0=1.
$$
Thus
$$
\int_Q|v_k|^2dV_0\leq
C.
$$
Since
$$
-\Delta_{g_k} v_k+c(n)R(g_k)v_k=c(n)R_k{ |u_k|}^\frac{4}{n-2}v_k,
$$
by Lemma \ref{epsilon},  $\|u_k\|_{C^0
(\epsilon,3L-\epsilon)\times S^{n-1}}\rightarrow 0$ for any $\epsilon>0$. By Lemma \ref{A1}, $v_k$ converges to a function $v$ in $W^{1,2}$, where $v$ solves the equation:
$$
-\Delta v+\frac{(n-2)^2}{4}v=0,
$$
and satisfies
$$
\int_{Q_2}|v|^2dV_0=1.
$$
Thus $v\neq 0$.

Moreover,
$$
\int_{[\epsilon,L]\times S^{n-1}}|v|^2dV_0=\lim_{k\rightarrow+\infty}
\int_{[\epsilon,L]\times S^{n-1}}|v_k|^2dV_0\leq e^{-L}
\lim_{k\rightarrow+\infty}\int_{Q_2}|v|^2dV_0,
$$
letting $\epsilon \rightarrow 0$ gives
$$
\int_{Q_1}|v|^2dV_0\leq e^{-L}\int_{Q_2}|v|^2dV_0.
$$
Similarly, there holds
$$
\int_{Q_3}|v|^2dV_0\leq e^{L}\int_{Q_2}|v|^2dV_0,
$$
which contradicts Lemma \ref{3-circle.R=0}. Hence, the statements in 1) are proved. Using the same arguments, we can easily carry out the proof of 2) and 3).
~\endproof

Next, we give two applications of Three Circles Theorem. The
first one is
the removability of singularity of scalar curvature equation.

\begin{cor}\label{singularity}
Let $g$ be a smooth metric over $B$, and $g'\in SC^p(B\setminus\{0\},g)$.
Assume
$$
\int_B(1+|R(g')|^p)dV_{g'}<+\infty
$$
for some $p>\frac{n}{2}$. Then  $g'$ can be extended to a metric in $SC^p(B,g)$.
\end{cor}

\proof Set $g'=u^\frac{4}{n-2}g$.
Choose a normal chart around $0$ with respect to $g$, and set
$$
(t,\theta)=(-\log r,\theta).
$$
On this polar coordinate, we set
$$
g'(t,\theta)=v^\frac{4}{n-2} g(t,\theta).
$$
Then
$$
v(t,\theta)=u(e^{-t},\theta)e^{-\frac{n-2}{2}t},\s
and\s \int_{B_{{\delta}}(0)}u^\frac{2n}{n-2}dx=
\int_{[-\log\delta,+\infty)\times S^{n-1}}v^\frac{2n}{n-2}dV_0.
$$
Without loss of generality, we assume
$$
\|g-g_0\|_{C^2((1,+\infty)\times S^{n-1})}<\tau,\s \int_{[0,\infty]\times S^{n-1}}
v^\frac{2n}{n-2}<\epsilon_0.
$$
By Lemma \ref{epsilon},  we may assume $\|v\|_{L^\infty((1,+\infty)\times S^{n-1})}<C$.
By { Three Circles Theorem, for any $t$, we have
$$
\int_{[t,t+1]\times S^{n-1}}v^2dV_0\leq Ce^{-\delta t}
\max\{\int_{[2t,2t+1]\times S^{n-1}}v^2dV_0,\int_{
[0,1]\times S^{n-1}}v^2dV_0\},
$$
where $\delta>0$ is a constant.
Then for any $p'\ge2$, we have
$$
\int_{[t,t+1]\times S^{n-1}}v^{p'}dV_0<C\int_{[t,t+1]\times
S^{n-1}}v^2dV_0<Ce^{-\delta t}.
$$
Then
$$
\int_{B_{e^{-t}}\setminus B_{e^{-(t+1)}}}u^{p'}dx
\leq\int_{[t,t+1]\times S^{n-1}}v^{p'}e^{\frac{p'(n-2)}{2}t-nt}dV_0
\leq C
e^{-\delta t-nt+(\frac{n-2}{2}p')t},
$$
i.e.
$$
\int_{B_{er}\setminus B_r}u^{p'}dx
\leq Cr^{
\delta+(n-\frac{n-2}{2}p')}.
$$
Choose $p'>\frac{2n}{n-2}$,
such that $\delta'=\delta+(n-\frac{n-2}{2}p')>0$.
Then
$$
\int_{B_{r_0}}u^{p'}dx=
\sum_k\int_{B_{2^{-k}r_0}\setminus B_{2^{-k-1}r_0}}u^{p'}dx\leq C\sum_k (2^{-k}r_0)^{\delta'}<+\infty.
$$
Put $u_r=r^{\frac{n-2}{2}}u(rx)$. We have
$$
\int_{B_e\setminus B_1}|\nabla u_r|^2dx=
\int_{B_{er}\setminus B_r}|\nabla u|^2dx,\s
\int_{B_{e}\setminus B_1}u_r^\frac{2n}{n-2}dx=
\int_{B_{er}\setminus B_r}u^\frac{2n}{n-2}dx.
$$
By Lemma \ref{epsilon},
$$
\int_{B_e\setminus B_1}|\nabla u_r|^2dx\leq C
\left(\int_{B_{e^2}\setminus B_\frac{1}{e}}u_r^\frac{2n}{n-2}dx\right)^\frac{n}{n-2}<Cr^{\delta}.
$$
Then
$$
\int_{B_{r_0}}|\nabla u|^2dx=
\sum_k\int_{B_{2^{-k}r_0}\setminus B_{2^{-k-1}r_0}}|\nabla u|^2dx\leq C\sum_k (2^{-k}r_0)^{\delta}<+\infty.
$$
Take a
cut-off function $\eta$ which is 0 on
$B_1$ and 1 on $B_{2}^c$ and set $\eta_r=
\eta(\frac{x}{r})$. We have
$$
\|\nabla\eta_ru\|_{L^2(B_{r_0})}\leq \|\nabla u\|_{L^2(B_{r_0})}
+\|\frac{u}{r}\|_{L^2(B_{2r})}.
$$
Since
$$
\int_{B_{2r}}(\frac{u}{r})^2dx\leq
(\int_{B_{2r}}u^\frac{2n}{n-2}dx)^\frac{n-2}{n}
(\int_{B_{2r}}r^{-n}dx)^\frac{2}{n}\leq
C(\int_{B_{2r}}u^\frac{2n}{n-2}dx)^\frac{n-2}{n}.
$$
We may choose $r_k\rightarrow 0$, such that
$\eta_{r_k}u$ converges weakly in $W^{1,2}$,
which implies that $u\in W^{1,2}$.

In a similar way, we can prove that $u$ satisfies
the scalar curvature equation weakly in $B$.
Using the arguments of the proof of Lemma \ref{epsilon},we get $u\in W^{2,p}(B)$.

By Lemma \ref{positive}, we get $u>0$.
\endproof

The following corollary will play an essential role in section 6:

\begin{cor}\label{positive.of.delta}
Let $g=u^\frac{4}{n-2}g_{\R^n}$, which satisfies
$$
\int_{\R^n}(1+|R(g)|^pu^\frac{2n}{n-2})dx<+\infty.
$$
Then
$$
\lim_{R\rightarrow+\infty}\int_{B_R}R(g)u^\frac{n+2}{n-2}dx>0.
$$
\end{cor}

\proof
Let $y=\frac{x}{|x|^2}$, and $g=v^\frac{4}{n-2}(y)
\sum\limits_idy^i\otimes dy^i$. By the above corollary,
$v\in W^{2,p}(\R^n)$  and
\begin{equation}\label{u.infty}
u(x)=v(\frac{x}{|x|^2})|x|^{2-n}.
\end{equation}
Hence there exists $C>1$, s.t.
$$
\frac{1}{C}<|x|^{n-2}u(x)<C.
$$
Let $u_k(x)=R_k^{n-2}u(R_kx)$,
where $R_k\rightarrow+\infty$. By \eqref{u.infty},
$u_k\rightarrow v(0)|x|^{2-n}$
in $C^0(B_2)$. Using the equation
$$
\Delta_{g_{R^{n}}}u_{k}+c(n)Ru_{k}^{\frac{n+2}{n-2}}R_{k}^{-2}=0,
$$
we get that
$u_k$ converges to $v(0)|x|^{2-n}$ weakly in $W^{2,p}
(B_2\setminus B_\frac{1}{2})$.
Hence by trace embedding, a subsequence of $\int_{\partial B_1}\frac{\partial u_k}{\partial r}dS$ converges to $(2-n)v(0)\omega_{n-1}$.

Since
$$
\int_{\partial B_1}\frac{\partial u_k}{\partial r}dS=
\int_{\partial B_{R_k}}\frac{\partial u}{\partial r}dS,
$$
Any subsequence of $\int_{\partial B_R}\frac{\partial u}
{\partial r}dS$ has a subsequence converges to $(2-n)v(0)\omega_{n-1}$.
Thus
$$
\lim_{R\rightarrow+\infty}\int_{\partial B_R}\frac{\partial u}{\partial r}dS
=(2-n)v(0)\omega_{n-1}.
$$
Then we get
$$
\int_{\R^n}Ru^\frac{n+2}{n-2}dx=\lim_{R\rightarrow+\infty}
\int_{B_R}Ru^\frac{n+2}{n-2}dx=-\lim_{R\rightarrow+\infty}\int_{
\partial B_R}\frac{\partial u}{\partial r}dS=(n-2)v(0)\omega_{n-1}.
$$
\endproof

\section{Bubble tree}
In this section, we set $(M,g)$ to be a smooth
closed Riemannian manifold, and
$g_k=u_k^\frac{4}{n-2}g\in SC^p_\Lambda(g)$ with $\vol(M,g_{k})= 1$.
Since
$$
\int_M|\nabla_gu_k|^2dV_g+c(n)\int_MR(g)u_k^2dV_g=c(n)\int_MR(g_k)u_k^\frac{2n}{n-2}dV_g\leq
c(n)\|R(g_k)\|_{L^1(M,g_k)},
$$
we have
$$
\int_M|\nabla_gu_k|^2dV_g\leq c(n)\left(\int_Mu_k^\frac{2n}{n-2}dV_g\right)^\frac{n-2}{n}
 \left(\int_M|R(g)|^\frac{n}{2}dV_g\right)^\frac{2}{n}+c(n)\|R(g_k)\|_{L^1(M,g_k)}.
$$
Therefore, $u_k$ is bounded in $W^{1,2}$ and
we may assume
$u_k$ converges weakly to $u_0$ in $W^{1,2}(M,g)$.

Throughout this section, we will denote the geodesic ball on $(M,g)$ of
radius $r$ centered at $p$ by $B_r(p)$.
First of all, we prove the following:

\begin{lem}
After passing to a subsequence, there exists  a finite
set $\mathcal{S}$, such that
$$
\lim_{r\rightarrow 0}\limsup_{k\rightarrow+\infty}
\vol(B_r(y),g_k)\leq\frac{\epsilon_0}{2},\s\forall
y\notin\CS,
$$
and
$$
\lim_{r\rightarrow 0}\liminf_{k\rightarrow+\infty}
\vol(B_r(y),g_k)>\frac{\epsilon_0}{2},\s\forall
y\in\CS.
$$
\end{lem}

\proof
Since $\int_Mu_k^\frac{2n}{n-2}dV_g=1$, by choosing a suitable subsequence, one  may assume
$u_k^\frac{2n}{n-2}dV_g$ converges to a Radon measure $\mu$
in the sense of distribution. Let
$$
\CS=\{y\in M:\mu(\{y\})>\frac{\epsilon_0}{2}\}.
$$
Obviously, $\CS$
is a finite set.

Take a  cut-off function $\eta$ which is 0 on $B_{r}^c(y)$, 1 on $B_\frac{r}{2}(y)$, and positive
on $B_r\setminus B_\frac{r}{2}(y)$. Then
$$
\vol(B_\frac{r}{2}(y),g_k)\leq\int_{M}\eta u_k^\frac{2n}{n-2}dV_g\leq \vol(B_{r}(y),g_k).
$$
Noting that $\int_{M}\eta u_k^\frac{2n}{n-2}dV_g\rightarrow\int_M\eta d\mu$, we complete the proof easily.
\endproof

Usually, we call $\CS$ the concentration set of $\{g_k\}$.
By Lemma \ref{epsilon}, $u_k$ converges weakly in
$W^{2,p}_{loc}(M\setminus\mathcal{S},g)$. Obviously, when  $\mathcal{S}=\emptyset$, $u_k$
converges to $u$ weakly in $W^{2,p}$ and  $u_0\in SC^p_\Lambda(M,g)$  with $\int_Mu_0^\frac{2n}{n-2}dV_g=1$. Hence, from now on,  we assume $\mathcal{S}\neq\emptyset$.

Let $x_0\in\mathcal{S}$ and
choose a normal chart $x_0$, such that $x_0=0$,
and
$g_{ij}=\delta_{ij}+O(|x|^2)$.

\subsection{The first bubble}
In this subsection, we show that there exists at least one bubble
at $x_0$.

Let $x_k\in B_\delta(0)$, $r_k>0$ such that
$$
\vol(B_{r_k}(x_k))=\frac{\epsilon_0}{2},\s
\vol(B_r(x))\leq \frac{\epsilon_0}{2},\s
\forall x\in B_\delta(0), r\leq r_k.
$$
Set $v_k=r_k^{\frac{n-2}{2}}u_k(x_k+r_kx)$, $\hat{g}_k(x)=g_{ij}(x_k+r_kx)dx^i\otimes dx^j$.
Then we have
$$
-\Delta_{\hat{g}_k}v_k+c(n)R(\hat{g}_k)v_{k}=c(n)R_{k}(x_k+r_kx)v_k^\frac{n+2}{n-2}.
$$
Obviously, we have
$$
\int_{B_r(x)}v_k^\frac{2n}{n-2}dV_{\hat{g}_k}=\int_{B_{rr_k}(rr_kx+x_k)}
u_k^\frac{2n}{n-2}dV_g,
$$
and
$$
\int_{B_r(x)}
|R_{k}|^p(x_k+r_kx)v_k^\frac{2n}{n-2}dV_{\hat g_k}=\int_{B_{rr_k}(rr_kx+x_k)}
|R_k(x)|^{p}u_k^\frac{2n}{n-2}dV_g.
$$
Then for any fixed $R>0$, we have
$$
\int_{B_1(x)}v_k^\frac{2n}{n-2}dx<\epsilon_0,\s \forall x\in B_R.
$$
Hence $v_k$ converges to a function $v$ weakly in $W_{loc}^{2,p}
(\R^n)$. We have
$$
\vol(\R^n,v^\frac{4}{n-2}g_{\R^n})\leq 1,\s \vol(B_1^{\R^n},v^\frac{4}{n-2}g_{\R^n})=\epsilon_0/2,\s \int_{\R^n}|R(v^\frac{4}{n-2}g_{\R^n})|^pv^\frac{2n}{n-2}\leq \Lambda.
$$
By Corollary \ref{singularity}, we can consider  $v^\frac{4}{n-2}g_{\R^n}$  as a metric in $SC^p(\Sp^n)$. We usually call $v$ or $v^\frac{4}{n-2}g_{\R^n}$
the first bubble.

It is easy to check that $\frac{1}{C}<{r_k}^{\frac{n-2}{2}}\|u_k\|_{L^\infty(B_\delta(x_0))}<C$ for some $C>0$.

\subsection{Some identities}
In this subsection, we assume that  for any $\epsilon>0$, we can
find $r$, such that
\begin{equation}\label{epsilon.neck}
\int_{B_{2t}(x_{k})\setminus B_t(x_k)}u_k^\frac{2n}{n-2}dx<\epsilon
\end{equation}
for any $t\in [\frac{r_k}{r},r]$.
By the arguments in \cite{Ding-Tian}, the above statement implies that
$u_k$ has only one bubble.
The main goal of
this subsection is to prove the following in this specical case, the
general case will be discussed in the next subsection.

\begin{lem}\label{identity.one bubble}
If \eqref{epsilon.neck} holds, then we have
\begin{equation}\label{volume.identity}
\lim_{r\rightarrow 0}\lim_{k\rightarrow+\infty}
\vol(B_r(x_k)\setminus B_\frac{r_k}{r}(x_k),g_k)=0,
\end{equation}
\begin{equation}\label{no.neck}
\lim_{r\rightarrow 0}\lim_{k\rightarrow+\infty}
diam_{g_k}(B_r(x_k)\setminus B_\frac{r_k}{r}(x_k),g_{k})=0,
\end{equation}
\begin{equation}\label{Yamabe.identity}
\lim_{r\rightarrow 0}\lim_{k\rightarrow+\infty}
\int_{B_r(x_k)\setminus B_\frac{r_k}{r}(x_k)}|R_k|u_k^\frac{2n}{n-2}dx=0.
\end{equation}
Moreover,  if $q\in (1,\frac{n}{n-2})$ and
$$
\int_{B_t(x)}|\phi_k|^qdx<Ct^{n-(n-2)q},
$$
holds for some  and any $B_t(x)$,
then
\begin{equation}\label{mk.identity}
\lim_{r\rightarrow 0}\lim_{k\rightarrow+\infty}
\int_{B_r(x_k)\setminus B_\frac{r_k}{r}(x_k)}|R_k\phi_k|u_k^\frac{4}{n-2}dx=0.
\end{equation}
\end{lem}

\proof
Set $T_k=-\log \frac{r_{k}}{r}$, $T=-\log r$ and
$$
\psi_k(t,\theta)=(-\log r,\theta) +x_k.
$$
On this polar coordinate, we set
$$
g_k(t,\theta)=v_k^{\frac{4}{n-2}}\psi_k^*(g)(t,\theta).
$$
Then
$$
v_k(t,\theta)=u_k((e^{-t},\theta)+x_k)e^{-\frac{n-2}{2}t}.
$$
We have
$$
\|\psi_k^*(g)-g_0\|_{C^2}\rightarrow 0,\s \int_{[t,t+3]\times S^{n-1}}
v_k^\frac{2n}{n-2}dx<\epsilon_0,\s whenever\s t\in [T,T_k-T],
$$
where $g_{0}$ is defined as in  section 4.
By Lemma \eqref{epsilon},  $\|v_k\|_{L^\infty}<C$.

 Without loss of generality, we assume
$T_{k}=T+m_k L$, where $m_k$ is a positive integer, and $L$ is a large fixed positive number as in section 4.
Putting $Q_i=[T+iL,T+(i+1)L]\times S^{n-1}$,
by Theorem \ref{3-circle}, we have
$$
\int_{Q_i}v_k^2dV_0\leq CA_k(L)(e^{-\delta iL}+e^{-\delta(m_k-i)L}),
$$
where $\delta>0$ and
$$
A_k(L)=\int_{[T,T+L]\times S^{n-1}}v_k^2dV_0+\int_{[T_k-L,T_k]\times S^{n-1}}
v_k^2dV_0.
$$
By \eqref{epsilon.neck}, we may choose $r$, such that
$$
A_k(L)<\epsilon'.
$$
Since
$$
\int_{Q_i}v_k^\frac{2n}{n-2}dx
<C\int_{Q_i}v_k^2dV_0<CA_k(L)(e^{-\delta i L}+e^{-\delta(m_k-i)L}),
$$
we get \eqref{volume.identity} by letting $\epsilon'\rightarrow 0$.\\

By elliptic estimate,
$$
diam(Q_i,g_k)<C\|v_k\|_{L^\infty(Q_i)}^{\frac{2}{n-2}}\leq \left(CA_k(L)(e^{-\delta i L}+e^{-\delta(m_k-i)L}\right))^\frac{1}{n},
$$
which implies \eqref{no.neck}.\\

From
$$
\int_{Q_i}|R_k|v_k^\frac{2n}{n-2}dx\leq
\left(\int_{Q_i}|R_k|^pv_k^\frac{2n}{n-2}dx\right)^\frac{1}{p}
\left(\int_{Q_i}v_k^\frac{2n}{n-2}dx\right)^\frac{p-1}{p}\leq
\Lambda(CA_k(L)(e^{-\delta i L}+e^{-\delta(m_k-i)L})^\frac{p-1}{p},
$$
it follows  \eqref{Yamabe.identity}.\\

Lastly, we prove \eqref{mk.identity}.
By
$$
\int_{B_{e^{-iL}r}\setminus B_{e^{-(i+1)L}r}(x_k)}|\phi_k|^qdx
\leq C{(e^{-iL}r)}^{n-(n-2)q},
$$
we have
$$
\int_{Q_i}|\phi_k|^qdx
\leq C(e^{-T-iL})^{n-(n-2)q}e^{n(T+iL)}=Ce^{-(T+iL)(n-2)q}.
$$
Then
\begin{eqnarray*}
\int_{B_{e^{-iL}r}\setminus B_{e^{-(i+1)L}r}(x_k)}|R_k\phi_k|u_k^\frac{4}{n-2}dx
&=&\int_{Q_i}|R_k\phi_k|(v_k)^\frac{4}{n-2}e^{-(n-2)t}dV_0\\
&=&
\int_{Q_i}|R_k|v_k^\frac{2n}{(n-2)p}
v_k^{\frac{4}{n-2}-\frac{2n}{(n-2)p}}|\phi_k|v_k e^{-(n-2)t}dV_0\\
&\leq&Ce^{-(n-2)(T+iL)}\left(\int_{Q_i}|R_k|^pv_k^\frac{2n}{n-2}dV_0
\right)^\frac{1}{p}\left(\int_{Q_i}|\phi_k|^qdV_0\right)^\frac{1}{q}\\
&&\times\left(\int_{Q_i}v_k^{\frac{4p-2n}{(n-2)p}s}dV_0\right)^\frac{1}{s}\\
&\leq&C\left(\int_{Q_i}|R_k|^pv_k^\frac{2n}{n-2}
dV_0\right)^\frac{1}{p}
\left(\int_{Q_i}v_k^{\frac{4p-2n}{(n-2)p}s}dV_0\right)^\frac{1}{s}\\
&\leq&C\Lambda
\|v_k\|_{L^\infty(Q_i)}^{\frac{4p-2n}{(n-2)p}}.
\end{eqnarray*}
We complete the proof.
\endproof

\subsection{Bubble Tree}
In this subsection, we will use the ideas in \cite{Chen-Li}
to construct bubble tree.

\begin{defi}
A sequence $\{(x_k,r_k): x_k\in B_\delta(0), r_k>0\}$ is called a  blowup sequence of $\{u_k\}$  if
$x_k\rightarrow 0$ and $r_k\rightarrow 0$, $v_k(x)=r_{k}^{\frac{n-2}{2}}u_k(x_k+r_kx) \rightharpoonup v\,\, \mbox{in}\,\, W^{2,p}_{loc}(\R^n\setminus S)$,
where $S$ is a finite set, and $v>0$.
\end{defi}

By removability  of singularity, $(\R^n,v^\frac{4}{n-2}g)$ can be extended to a metric in  $SC^p(\mathbb{S}^n,g_{\mathbb{S}^n})$.

\begin{defi}\label{essen}
Two blowup sequences $\{(x_k,r_k)\}$ and $\{(x_k',r_k')\}$ of $\{u_k\}$  are said to be
{\it essentially different} if one of the following happens
\begin{equation}\label{bubble1}
\frac{r_k}{r_k'}\rightarrow+\infty,\mbox{ or }
\frac{r_k'}{r_k}\rightarrow+\infty,
\s \mbox{or} \s\frac{|x_k-x_k'|}{r_k+r_k'}
\rightarrow+\infty.
\end{equation}
Otherwise, they are called  essentially same.
\end{defi}
In the sequel, we will write $(x_k,r_k)$ for a blowup sequence. For simplicity, we set $(0,1)$ to be the $0$-blow up sequence.

Let $v_k^\alpha(x)=(r_{k}^{\alpha})^{\frac{n-2}{2}}u_k(x_k^\alpha+r_k^\alpha x)$. We assume
$v_k^\alpha$ converges to $v^\alpha$ in $W^{2,p}(\R^n\setminus\CS^\alpha)$ weakly,
where $\CS^\alpha$ is the set of concentration points of $\{v_k^\alpha\}$, which is finite. Set
$$
U^\alpha_r=B_\frac{1}{r}\backslash \bigcup_{p\in \CS^\alpha}B_r(p),\s \Omega^\alpha_{r,k}=x_k^\alpha+r_k^\alpha U^\alpha_r.
$$
It is easy to check that, given $m$-bubbles, after passing
to a subsequence, we have
$$
\Omega^\alpha_{r,k}\cap\Omega^\beta_{r,k}=\emptyset,
$$
when $r$ is sufficiently small and $k$ is sufficiently large. For more details, one can refer to \cite{Chen-Li}.
Since
$$
\lim_{r\rightarrow 0}\lim_{k\rightarrow+\infty}\int_{U_r^\alpha}(v_k^\alpha)^\frac{2n}{n-2}=\int_{\R^n}(v^\alpha)^\frac{2n}{n-2}dx\geq \tau,
$$
we have
$$
m\tau\leq\Lambda^\frac{1}{p}.
$$
Thus, there exists only finitely many bubbles.
Usually, we say the sequence  $\{u_k\}$ has $m$ bubbles if $\{u_k\}$ has $m$ essentially different blowup sequences and no subsequence of $\{u_k\}$ has more than $m$ essentially different blowup sequences.

Now, we assume $\{u_k\}$ has $m$ essentially different blowup sequences $(x_k^1,r_k^1)$,
$\cdots$, $(x_k^m,r_k^m)$  at $0$. It is convenient to assume that:
$\frac{r_k^\alpha}{r_k^\beta}\rightarrow 0$ for any $\alpha\not=\beta$, we always assume, up to selecting a subsequence if needed, that
\begin{equation}\label{bubble2}
\mbox{either}\s\frac{|x_k^\alpha-x_k^\beta|}{r_k^\beta}\rightarrow+\infty,\s
\mbox{or}\s\frac{x_k^\alpha-x_k^\beta}{r_k^\beta}\s \mbox{converges.}
\end{equation}

\begin{defi}\label{less}
For two essentially different blowup sequences $\{(x^\alpha_k,r^\alpha_k)\}, \{(x^\beta_k,r^\beta_k)\}$, we say $(x_k^\alpha,r_k^\alpha)<(x_k^\beta,r_k^\beta)$, if
$\frac{r_k^\alpha}{r_k^\beta}\rightarrow 0$ and
$\frac{x_k^\alpha-x_k^\beta}{r_k^\beta}$ converges as $k\to\infty$.
\end{defi}

By arguments in \cite{Chen-Li},
\begin{equation}\label{equiv}
(x_k^\alpha,r_k^\alpha)<(x_k^\beta,r_k^\beta)
\Longleftrightarrow D_{Rr_k^\alpha}(x_k^\alpha)\subset D_{Rr_k^\beta}(x_k^\beta),\s\mbox{for some $R$ and all large $k$.}
\end{equation}
and the bubble $v^\alpha$ of $\{u_k\}$ at 0 is also a bubble of $\{v_k^\beta\}$
at $x_{\alpha\beta}= \lim_{k\rightarrow+\infty}\frac{x_k^\alpha-x_k^\beta}{r_k^\beta}$.

\begin{defi}\label{ontop}
A blowup sequence $(x_k^\alpha,r_k^\alpha)$ is said to be {\it right on top of} another blowup sequence $(x_k^\beta,r_k^\beta)$, if
$(x_k^\alpha,r_k^\alpha)<(x_k^\beta,r_k^\beta)$ and there is no blowup sequence
$(x_k^\gamma,r_k^\gamma)$ which is essentially different from $(x^\alpha_k,r^\alpha_k)$ and $(x^\beta_k,r^\beta_k)$, such that
$(x_k^\alpha,r_k^\alpha)<(x_k^\gamma,r_k^\gamma)<(x_k^\beta,r_k^\beta).$
We say $(x_k^\alpha,r_k^\alpha)$ is a treetop if there is no
blowup sequence on the top of it.
\end{defi}

We write $B_r(x_0)\setminus (\cup_\alpha\Omega_{r,k}^\alpha)$
as the  disjoint union of its components:
$$
B_r(x_0)\setminus (\bigcup_\alpha\Omega_{r,k}^\alpha)=\bigcup_\beta N_{r,k}^\beta.
$$
It is easy check that $N_{r,k}^\beta$ is topologically
a ball minus finitely many small balls (see \cite{Chen-Li} for more details).

\begin{pro}
We have
\begin{equation*}
\lim_{r\rightarrow 0}\lim_{k\rightarrow+\infty}
vol(N_{r,k}^\beta,g_k)=0,
\end{equation*}
\begin{equation*}
\lim_{r\rightarrow 0}\lim_{k\rightarrow+\infty}
diam_{g_k} N_{r,k}^\beta=0,
\end{equation*}
\begin{equation}\label{Yamabe.identity2}
\lim_{r\rightarrow 0}\lim_{k\rightarrow+\infty}
\int_{N_{r,k}^\beta}|R_k|u_k^\frac{2n}{n-2}=0.
\end{equation}
Moreover,  if $q\in (1,\frac{n}{n-2})$ and
\begin{equation}\label{Lq.in.Br}
\int_{B_t(x)}|\phi_k|^qdx<Ct^{n-(n-2)q}
\end{equation}
for any $B_t(x)$,
then
\begin{equation}\label{mk.identity2}
\lim_{r\rightarrow 0}\lim_{k\rightarrow+\infty}
\int_{N_{r,k}^\beta}|R_k\phi_k|u_k^\frac{4}{n-2}dx=0.
\end{equation}
\end{pro}

\proof We only prove \eqref{mk.identity2}.
We will prove it by induction of $m$, the number of bubbles.\\
when $m=1$, we have only one bubble. If (5.1) is not true, then their
exists $\epsilon'$ such that $\forall r$, $\exists$
$t_{k}\in[\frac{r_{k}}{r},r]$,
\begin{equation}
\int_{B_{2t_{k}}(x_{k})\setminus B_{t_{k}}(x_k)}u_k^\frac{2n}{n-2}dx=\epsilon'.
\end{equation}
As $r\rightarrow 0$, $\frac{t_{k}}{r_{k}}\rightarrow \infty$. Then
either $(x_k,t_k)$ is  another bubble or $r_{k}^{\frac{n-2}{2}}u_k(r_kx+x_k)$
has at least one concentration point. For the latter case, $r_{k}^{\frac{n-2}{2}}u_k(r_kx+x_k)$
has at least one bubble.
Altogether,  $u_k$ has at least 2  bubbles, a contradiction. Then by lemma 5.2 we have the result.

We assume the result holds for the case of $m-1$ bubbles, we prove the case of $m$ bubbles. Consider the following two cases separately. Case 1: we have only one essentially different blow up sequence, denoted by $(x_{k}^{1},r_{k}^{1})$, which is right top of $(0,1)$; Case 2:  there are  essentially different blow up sequences $(x_{k}^{\alpha},r_{k}^{\alpha})$, $\alpha=1,2,...m'$, right on top of $(0,1)$, where $m'>1$.

For Case 1,
other $m-1$ bubbles can be considered as   bubbles of
$$
u_{k}^{1}=(r_{k}^{1})^{\frac{n-2}{2}}u_{k}(r_{k}^{1}x+x_{k}^{1}).
$$
It is easy to check that
$\lim\limits_{k\rightarrow \infty}\dint_{B_{2t}(x_{k}^{1})\setminus B_{t}(x_{k}^{1})}u_{k}^{\frac{2n}{n-2}}dx<\epsilon$ for all $t\in [\frac{r_{k}^{1}}{r},r]$ when $r$ is sufficiently small.
Then using the argument in the proof of Lemma \ref{identity.one bubble},
we can prove
$$
\lim_{r\rightarrow 0}\lim_{k\rightarrow+\infty}
\int_{ B_r(x_k^1)\setminus B_{\frac{r_k^1}{r}}(x_k^1)}|R_k\phi_k|u_k^\frac{4}{n-2}dx=0.
$$
Then by induction, we get the result.\\

For Case 2,
it is easy to check that $\frac{|x_{k}^{\alpha}-x_{k}^{\beta}|}{r_{k}^{\alpha}} \rightarrow +\infty$  if $\alpha \neq \beta$, and $\alpha,\beta\leq m'$.
Without loss of generality , we assume that \\
$$
r_{k}'=\max \{|x_{k}^{\alpha}-x_{k}^{\beta}|:\alpha,\beta =1,...m'\}=|x_{k}^{1}-x_{k}^{2}|,
$$
then $\frac{r_{k}'}{r_{k}^{\alpha}}\geq\max \{\frac{|x_{k}^{1}-x_{k}^{\alpha}|}{r_{k}^{\alpha}},\frac{|x_{k}^{2}-x_{k}^{\alpha}|}{r_{k}^{\alpha}}\} \rightarrow +\infty$, as $k\rightarrow +\infty$. Every bubble $v^{\alpha}$ is a bubble of $u_{k}'={r_{k}'}^{\frac{n-2}{2}}u_{k}(r_{k}'x+x_{k}^{1})$ because
\begin{eqnarray*}
u_{k}^{\alpha}(x)&=&{r_{k}^{\alpha}}^{\frac{n-2}{2}}{ u_{k}}(r_{k}^{\alpha}x+x_{k}^{\alpha})\\
&=&(\frac{r_{k}^{\alpha}}{r_{k}'})^{\frac{n-2}{2}}{r_{k}'}^{\frac{n-2}{2}}
u_{k}(x_{k}^{1}+r_{k}'(\frac{x_{k}^{\alpha}-x_{k}^{1}}{r_{k}'}+\frac{r_{k}^{\alpha}}{r_{k}'}x))\\
&=&(\frac{r_{k}^{\alpha}}{r_{k}'})^{\frac{n-2}{2}}u_{k}'(\frac{x_{k}^{\alpha}-x_{k}^{1}}{r_{k}'}+\frac{r_{k}^{\alpha}}{r_{k}'}x).
\end{eqnarray*}
However, $(x_k^1,r_k')$ is not blowup sequence,
$u_{k}'$ must converge to 0 in $W^{1,2}$ weakly.
Moreover, it is easy to check that
$\lim\limits_{k\rightarrow \infty}\dint_{B_{2t}(x_{k}^{1})\setminus B_{t}(x_{k}^{1})}u_{k}^{\frac{2n}{n-2}}dx<\epsilon$ for all $t\in [\frac{r_{k}^{1}}{r},r]$ when $r$ is sufficiently small.

 By the choice of $r_{k}'$, $u_{k}'$ has at least 2 blow up points.
By taking two small enough balls around the two blow up points, at each of them, there are at most $m-1$ bubbles. Using the proof for case $m=1$ and  induction to the two blow up points separately, we finish the proof.
\endproof

\section{The proof of Theorem \ref{main2}}
\noindent{\it Proof of Theorem \ref{main2}.}  For simplicity, we assume $\mathcal{S}$
contains only one point $x_0$.

We set
\begin{equation}\label{def.ck}
\frac{1}{c_{k}}=\int_{M}|c(n)R_{k}u_{k}^{\frac{n+2}{n-2}}-c(n)R(g)u_{k}+u_{k}|dV_{g}
\end{equation}
Noting that
$$
-\Delta_{g} u_{k} +u_{k}=c(n)R_{k}u_{k}^{\frac{n+2}{n-2}}-(c(n)R-1)u_{k},
$$
by Lemma \ref{Lq}, we get
$$
\int_{M}|c_{k}u_{k}|^{q}dV_g<C(q,M,g),\s \int_{M}|\nabla_g c_{k}u_{k}|^{q'}dV_g<C(q',M,g),
$$
where $q\in(1,\frac{n}{n-2})$ and $q'\in(1,\frac{n}{n-1})$. Then we may assume $c_{k}u_{k}$ converges weakly in $W^{1,q'}$ to a function $G$. By Lemma \ref{A1} , we also know $c_{k}u_{k}$ converges to $G$ in $W^{2,p}_{loc}(M\setminus\CS)$ weakly.

During the rest of this section, we will study the behavior of $c_ku_k$ near
$x_0$.
We have two possibilities: i) $\lim\limits_{k\rightarrow \infty}\dint_{M}c_{k}|R_{k}|u_{k}^{\frac{n+2}{n-2}}dV_g=0$; ii) $\lim\limits_{k\rightarrow \infty}\dint_{M}c_{k}|R_{k}|u_{k}^{\frac{n+2}{n-2}}dV_g>0$.

For case i):
we have
$$
\int_{M}|G-c(n)RG| =1.
$$
In distribution sense, $G$ satisfies
$$
-\Delta_{g} G+c(n)R(g)G=0.
$$
Then $G$ is smooth by elliptic equation estimate and the scalar curvature of $(M,G^{\frac{4}{n-2}}g)$ is $0$.

For case ii):
We claim that $G$ is  some kind of Green function
which satisfies
$$
-\Delta_{g}F+c(n)R(g)F = \lambda \delta_{x_0},\s \lambda>0.
$$
Clearly, we only need to prove
$$
\lim\limits_{k\rightarrow\infty}\dint_{M}c(n)\varphi R_{k}u_{k}^{\frac{n+2}{n-2}}c_{k}dV_g=\lambda \varphi(p),\s \forall \varphi\in C^{\infty}(M).
$$

We divide the integral $\dint_{M}c(n)\varphi R_{k}u_{k}^{\frac{n+2}{n-2}}c_{k}dV_g$ into three parts:
$$
\int_{M\setminus B_r(x_0)} \varphi R_{k}u_{k}^{\frac{n+2}{n-2}}c_{k}dV_g+\int_{\bigcup_{\beta} N_{r,k}^\beta} \varphi R_{k}u_{k}^{\frac{n+2}{n-2}}c_{k}dV_g+\int_{\bigcup_\alpha\Omega_{r,k}^\alpha}\varphi R_{k}u_{k}^{\frac{n+2}{n-2}}c_{k}dV_g
$$
where $N_{r,k}^{\beta}$ and $\Omega_{r,k}^\alpha$ are as in section 5.

Since  $c_ku_{k}$ converges to $G$ and $u_k$  converges to $0$ uniformly on $M\setminus B_{r}(x_0)$,
 we  have
{\begin{eqnarray*}
&&\lim_{k\rightarrow \infty}\int_{M\setminus B_r(x_0)}|R_{k}|u_{k}^{\frac{n+2}{n-2}}c_{k}dV_g\\
&=&{\lim_{k\rightarrow \infty}\int_{M\setminus B_r(x_0)}|R_{k}|u_{k}^{\frac{2n}{n-2}\frac{q-1}{q}}u_{k}^{\frac{4}{n-2}-\frac{2n}{n-2}\frac{q-1}{q}}(c_{k}u_{k})}dV_g\\
&\leq&\lim_{k\rightarrow \infty}
 C(M,q,g)^{\frac{1}{q}}\left(\int_{M\setminus B_r(x_0)}|R_k|^\frac{q}{q-1}dV_{g_k}\right)^\frac{q-1}{q} \max\limits_{M\setminus B_{r}(x_{0})}u_{k}^{\frac{4}{n-2}-\frac{2n}{n-2}\frac{q-1}{q}}\\
 &\leq&\lim_{k\rightarrow \infty}C(M,q,g)^{\frac{1}{q}}\left(\int_{M\setminus B_r(x_0)}|R_k|^pdV_{g_k}\right)^\frac{1}{p}
\left(\int_{M\setminus B_r(x_0)}u_k^\frac{2n}{n-2}dV_g\right)^{\frac{q-1}{q}-\frac{1}{p}} \max\limits_{M\setminus B_{r}(x_{0})}u_{k}^{\frac{4}{n-2}-\frac{2n}{n-2}\frac{q-1}{q}}\\
&\leq&\lim_{k\rightarrow \infty}C(M,q,g)^{\frac{1}{q}}\Lambda^{\frac{1}{p}}
\left(\int_{M\setminus B_r(x_0)}u_k^\frac{2n}{n-2}dV_g\right)^{\frac{q-1}{q}-\frac{1}{p}}\max\limits_{M\setminus B_{r}(x_{0})}u_{k}^{\frac{4}{n-2}-\frac{2n}{n-2}\frac{q-1}{q}}\\
&=&0.
\end{eqnarray*}

where $\frac{p}{p-1}<q<\frac{n}{n-2}$.

Since $c_ku_k$ satisfies \eqref{Lq.in.Br},
by \eqref{mk.identity2}, we have
$$
\lim_{r\rightarrow 0}\lim_{k\rightarrow+\infty}\int_{\bigcup_{\beta}N_{k,r}^\beta}
|R_k|u_k^\frac{n+2}{n-2}c_kdx=0,
$$
which yields that
\begin{equation}\label{|Rk|}
\lim_{k\rightarrow+\infty}\int_M|R_k|u_k^\frac{n+2}{n-2}c_kdx=\lim_{r\rightarrow0}\lim_{k\rightarrow+\infty}\int_{\bigcup_\alpha\Omega_{r,k}^\alpha}|R_{k}|u_{k}^{\frac{n+2}{n-2}}c_{k}dx,
\end{equation}
and
$$
\lim_{r\rightarrow 0}\lim_{k\rightarrow+\infty}\int_{\bigcup_{\beta}N_{r,k}^{\beta}} c(n)R_ku_k^\frac{n+2}{n-2}c_k\varphi dx=0.
$$

Next, we estimate
$\lim\limits_{r\rightarrow0}\lim\limits_{k\rightarrow+\infty}\int_{\Omega_{r,k}^\alpha}\varphi R_{k}u_{k}^{\frac{n+2}{n-2}}c_{k}
$ step by step.\\

\noindent{\bf Claim:} $
\sum\limits_\alpha\lim_{k\rightarrow+\infty}(c_k(r_k^\alpha)^\frac{n-2}{2})>0.
$\\

Since $\int_M(c_ku_k)^{q}<C$, it follows from \eqref{def.ck} that
$$
0<\lim_{r\rightarrow0}\lim_{k\rightarrow+\infty}\int_{\cup_\alpha\Omega_{r,k}^\alpha}|R_{k}|u_{k}^{\frac{n+2}{n-2}}c_{k}dx<C.
$$

Let
\begin{eqnarray*}
\lambda^{\alpha}&=&\lim_{r\rightarrow 0}\lim_{k\rightarrow+\infty}\int_{\Omega_{r,k}^\alpha}
|R_k|u_k^\frac{n+2}{n-2}c_kdx\\
&=&\lim_{r\rightarrow 0}\lim_{k\rightarrow+\infty}\int_{U_{\alpha}^r}|R_k(r_k^\alpha x+x_k^\alpha)|(u_k^\alpha)^\frac{n+2}{n-2}c_k(r_k^\alpha)^\frac{n-2}{2}dx,
\end{eqnarray*}
 where $U^{\alpha}_{r}=D_{\frac{1}{r}}(0)\setminus\bigcup\limits_{p\in\CS^\alpha}D_{r}(p)$ and
$\CS^\alpha$ is the concentration point set of $u_k^\alpha$. Set
$$
f_k^{\alpha+}=R_k^+(r_k^\alpha x+x_k^\alpha)(u_k^{\alpha})^\frac{n+2}{n-2},\s
f_k^{\alpha-}=R_k^-(r_k^\alpha x+x_k^\alpha) (u_k^\alpha)^\frac{n+2}{n-2}.
$$
We may assume $f_k^{\alpha+}$ and $f_k^{\alpha-}$ converges weakly to
$f^{\alpha+}$ and $f^{\alpha-}$ respectively in
$L^p(U_r^\alpha)$ for any $r$. Hence
$$
\lambda^\alpha=\lim_{r\rightarrow 0}\lim_{k\rightarrow+\infty}(c_k(r_k^\alpha)^\frac{n-2}{2})\int_{U_r^\alpha}(f^{\alpha+}+f^{\alpha-})dx.
$$
Since $\sum\lambda^\alpha>0$, we prove the claim.\\

\noindent{\bf Claim: }$
\sum\limits_\alpha\lim_{k\rightarrow+\infty}(c_k(r_k^\alpha)^\frac{n-2}{2})<+\infty.$\\

Let $R^\alpha$ is the scalar curvature of $(v^{\alpha})^\frac{4}{n-2}g_{\R^n}$, where $g_{\R^n}$ is standard Euclidean metric. Obviously,
$$
R^\alpha (v^{\alpha})^\frac{n+2}{n-2}=f^{\alpha+}-f^{\alpha-}.
$$
Noting that   Corollary \ref{singularity} implies that  $f^{\alpha+}-f^{\alpha-}\in L^p_{loc}(\R^n)$.

By Corollary \ref{positive.of.delta},
$$
\lim_{r\rightarrow 0}\int_{U_r^\alpha}(f^{\alpha+}-f^{\alpha-})dx>0,
$$
then
$$
\lim_{r\rightarrow 0}\int_{U_r^\alpha}(f^{\alpha+}+f^{\alpha-})dx>0,
$$
hence
$$
\sum_{\alpha}\lim_{k\rightarrow+\infty}c_k(r_k^\alpha)^\frac{n-2}{2}<+\infty.
$$

\noindent{\bf Claim: }$\lim\limits_{r\rightarrow0}\lim\limits_{k\rightarrow+\infty}\int_{\Omega_{r,k}^\alpha}\varphi R_{k}u_{k}^{\frac{n+2}{n-2}}c_{k}dx
=\lambda\varphi(x_0)$ for some $\lambda>0$.

We have
\begin{eqnarray*}
&&\lim_{r\rightarrow 0}\lim_{k\rightarrow+\infty}
\int_{\Omega_{r,k}^\alpha}\varphi R_ku_k^\frac{n+2}{n-2}c_kdx\\
&=&\sum_\alpha
\lim_{r\rightarrow0}\lim_{k\rightarrow+\infty}\int_{U_r^\alpha}\varphi(x_0) R_k(r_k^\alpha x+x_k^\alpha)(u_k^\alpha)^\frac{n+2}{n-2}(r_k^{\alpha})^{{\frac{n-2}{2}}}c_kdx\\
&=&\varphi(x_0)\sum_\alpha\lim_{k\rightarrow+\infty}(r_k^{\alpha})^{{\frac{n-2}{2}}}c_k\lim_{r\rightarrow 0}
\int_{U_r^\alpha} R^\alpha(v^\alpha)^\frac{n+2}{n-2}dx\\
&:=&\lambda\varphi(x_0)
\end{eqnarray*}
By Corollary \ref{positive.of.delta} again, $\lambda>0$.\\

\noindent{\bf Claim: } When we have only one bubble at $x_0$, we have $
\frac{1}{C}<\frac{c_k}{m_k}<C.
$ for some $C>0$.\\

If $m_k$ is bounded, by the equation of $u_k$, $\|u_k\|_{W^{2,p}}$
is bounded, then $\CS=\emptyset$. Hence $m_k\rightarrow+\infty$. Put $m_k=u_k(x_k)$. Obviously, $x_k\rightarrow x_0$.

Let $t_k^\frac{2-n}{2}=m_k$, $v_k(x)=u_k(x_k+t_kx)/m_k$. We have the equation
$$
-\Delta_{g(x_k+t_kx)}v_k+{ c(n)t_{k}^{2}}R(x_k+t_kx)v_k=c(n){R_{k}}(x_k+t_kx)v_k^\frac{n+2}{n-2},
$$
and
$$
\int_{B_R}|R(x_k+t_kx)|^pv_k^\frac{2n}{n-2}\leq\Lambda,\s  v_k(0)=\|v_k\|_{L^\infty}= 1.
$$
Then $v_k$ converges weakly to a function $v$ with $v(0)=1$. Therefore,
$(x_k,t_k)$ is a blowup sequence. Then we can find $r_k^\beta$, such that
$0<\frac{1}{C}<\frac{t_k}{r_k^\beta}<C$.  Since we have only one bubble, then
$$
0<\frac{1}{C}<t_k^\frac{n-2}{2}c_k<C,
$$
which implies the claim.

\section{The proof of Theorem \ref{main3}}

\noindent{\it Proof of Theorem \ref{main3}.} We will divide the proof
into three cases.\\

\noindent{\bf Case 1.}   $G$ is smooth at a concentration point
$x_0$.

Take a normal chart at $x_0$.
Let $(x_k,r_k)$
be a treetop,
and set
$v_k(x)=r_k^\frac{n-2}{2}u_k(r_kx+x_k)$ and assume $v_k$
converges to $v$. Let $g'=v^\frac{4}{n-2}g_{\R^n}$ and $g_0=G^\frac{4}{n-2}g$.

We define
$$\Gamma_t=
\{t\theta:\theta\in B_\rho^{\mathbb{S}^{n-1}}(\theta_0)\},
$$
where $B_\rho^{\mathbb{S}^{n-1}}(\theta_0)$ is the geodesic ball of
$\mathbb{S}^{n-1}$ of radius $r$ centered at $\theta_0$.

Obviously, we can find a ball $\rho$ and
$\theta_0$, such that
$$
\vol(g_k,\bigcup_{t\in [0,\delta]}\Gamma_t)<\epsilon_0.
$$
Take $t_{k}\in[\frac{r_k}{\delta},\delta]$, such that
\begin{equation*}
diam_{g_k}\Gamma_{t_k}
=\inf_{t\in[\frac{r_k}{\delta},\delta]}diam_{g_k}\Gamma_t.
\end{equation*}
Note that
$$
\lim_{k\rightarrow +\infty}diam_{g_k}(\Gamma_{\frac{r_k}{r}})=diam_{g'}({ \Gamma_{\frac{1}{r}}})>0,\s
\lim_{k\rightarrow+\infty}diam_{c_k^\frac{4}{n-2}g_k}(\Gamma_{r})
=diam_{g_0}(\Gamma_r)>0.
$$
By removability of singularity, we have
$$
\lim_{t\rightarrow +\infty}\lim_{k\rightarrow +\infty}diam_{g_k}(\Gamma_{tr_k})=0,\s
\lim_{t\rightarrow 0}\lim_{k\rightarrow+\infty}diam_{g_k}(\Gamma_t)=0.$$
Then
$$t_k\rightarrow 0,\s and \s \frac{t_k}{r_k}\rightarrow+\infty.$$

We set
$$
\tilde{v}_k=t_{k}^{\frac{n-2}{2}}u_{k}(x_{k}+t_{k}x), \s
\tilde{g}_k=\tilde{v}_{k}^\frac{4}{n-2}g_k(x_k+t_kx).
$$
Then
\begin{equation}\label{tk}
diam_{\tilde{g}_k}\Gamma_1
=\inf_{t\in[\frac{r_k}{t_{k}\delta},\frac{\delta}{t_k}]}diam_{\tilde{g}_k}\Gamma_t.
\end{equation}

Let $\tilde{\CS}$ be the concentration point set of $\tilde{v}_k$.
Obviously, $0\in \tilde{\CS}$. By th choice of $\Gamma_t$,
\begin{equation}\label{concentrationonGamma}
\tilde{\CS}\cap \bigcup_{t>0}\Gamma_t=\emptyset.
\end{equation}

Since $\lim_{k\rightarrow +\infty}diam_{\tilde{g}_k}(\Gamma_{1})\rightarrow
0$,
 $\tilde{v}_k$ converges to 0 weakly in  $W^{1,2}_{loc}(R^{n}\setminus \tilde{\CS})$.

Next, we shall show that there exists $\tilde{c}_k\rightarrow+\infty$,  such that
$\tilde{c}_k\tilde{v}_k(x)$ converges weakly in $W^{2,p}_{loc}(\R^n\setminus \hat{S})$ to a positive harmonic function. Before we continue the proof,
we need to define a radius concerning John-Nirenberg inequality.

Define
$$
\rho(u,x)=\sup \{r:t^{2-n}\int_{B_{t}(x)}|\nabla\log u|^2dx<\hat{\epsilon}^2_{0}, \s\forall  t<r\},
$$
where $\hat{\epsilon}_0>0$ is a fixed constant.
It is easy to check that when $u\in W^{2,p}$, $\rho(u,x)>0$.
The key
observation about $\rho(\tilde{v}_k,x)$ is the following:

\begin{lem}\label{brk}
For any $\Omega\subset\subset\R^n\setminus(\{0\}\cup\hat{\CS})$, there
exists a constant $a>0$, such that
\begin{equation*}
\inf_{x\in\overline{\Omega}} \rho(\tilde{v}_k,x)>a,
\end{equation*}
when $k$ is sufficiently large.
\end{lem}

\proof We denote $\rho(\tilde{v}_k,x)$ by $\rho_k(x)$.

Assume the lemma is not true. Then there exists $x_k'\in\overline{\Omega}$, such that $x_k'\rightarrow y_0$ and
$\rho_k(x_k')\rightarrow 0$.

Take  $\delta$, such that $B^{\R^n}_\delta(y_0)\cap\hat{\CS}=\emptyset$.
Put $y_k\in B^{\R^n}_\delta(y_0)$, such that
$$
\frac{\rho_k(y_k)}{\delta-|y_k-y_0|}=\inf_{y\in B_\delta(y_0)}\frac{\rho_k(x)}{\delta-|x-y_0|}:=\lambda_k.
$$
Note that $\lambda_k\leq \frac{\rho_k(y_k)}{\delta-|y_k-y_0|}\rightarrow 0$. We get $
\rho_k(y_k)=\lambda_k(\delta-|y_k-y_0|)\rightarrow 0$, and
$R\rho_k(y_k)\leq \delta-|y_k-y_0|$ for any fixed $R$ and sufficiently large $k$. Then
\begin{equation}\label{BR}
B_{R\rho_k}(y_k)\subset B_\delta(y_0),
\end{equation}
where $k$ is sufficiently large and $R$ is  fixed.
Moreover, for any $y\in B_{Rr_k}(y_k)$, we have
\begin{eqnarray*}
\frac{\rho_k(y)}{\rho_k(y_k)}&\geq& \frac{\delta-|y-y_0|}{\delta-|y_k-y_0|}\\
&\geq& 1+\frac{|y_k-y_0|-|y-y_0|}{\delta-|y_k-y_0|}\\
&\geq&
1-\frac{\left||y_k-y_0|-|y-y_0|\right|}{\delta-|y_k-y_0|}\\
&\geq& 1-\frac{|y-y_k|}{\delta-|y_k-y_0|}\\
&\geq&
1-\frac{Rr_k(y_k)}{\delta-|y_k-y_0|}.
\end{eqnarray*}
Hence
\begin{equation}\label{boundnessofradius}\frac{\rho_k(y)}{\rho_{k}(y_{k})}>\frac{1}{2}
\end{equation} when $k$ is sufficiently large.

Let $v_k'(x)=(\rho_k(y_k))^{\frac{n}{2}-1}\tilde{v}_k(y_k+\rho_k(y_k)x)$.
We have
$$
\int_{B_1}|\nabla \log v_k'|^2dx=\hat{\epsilon}^2_0.
$$
By \eqref{BR} and \eqref{boundnessofradius},
$\rho(v_k',x)\geq\frac{1}{2}$ for any $x\in B^{\R^n}_R(0)$ and sufficiently
large $k$.
Choose $c_k'$, such that $\int_{B_1}\log c_k'v_k'dx=0$. We can find
$C(R)$, such that $\int_{B_R}|\nabla\log v_k'|^2dx
<C(R)$. By Poincare inequality,
$\|\log c_k'v_k'\|_{W^{1,2}(B_R)}<C$. Then we may assume $\log c_k'v_k'$
converges to a function $\phi$ in $L^q(B_R)$ for any $q<\frac{2n}{n-2}$.  After choosing a subsequence still denoted by $\log c_k'v_k'$ ,we can assume $\log c_k'v_k'$ converges to $\phi$  pointwise   for a.e. $x\in\R^n$.
By Eropob Theorem, we can find a measurable set $F\subset B_{R}$, such that $L(B_R\setminus F)<\epsilon$ and
$\log c_{k}'v_{k}'$ converges to $\phi$ uniformly on $F$, where $L$ is the
Lebesgue measure.

By Theorem \ref{J-N}, we have
$$
\int_{B_{1}}e^{\frac{\alpha}{\hat{\epsilon}_0}\log c_{k}'v_{k}'}dx\int_{B_{1}}e^{-\frac{\alpha}{\hat{\epsilon}_0}\log c_{k}'v_{k}'}dx<\beta^{2}.
$$
Since $e^{\frac{\alpha}{\hat{\epsilon}_0}\log c_k'v_k'}$ and
$e^{-\frac{\alpha}{\hat{\epsilon}_0}\log c_k'v_k'}$   converge uniformly  to $e^{\frac{\alpha}{\hat{\epsilon}_0}\phi}$ and
$e^{-\frac{\alpha}{\hat{\epsilon}_0}\phi}$ respectively on $F$,
$$
\int_{B_R}e^{-\frac{\alpha}{\hat{\epsilon}_0}\log c_k'v_k'}dx\geq a_1>0,\s
and\s
\int_{B_R}e^{\frac{\alpha}{\hat{\epsilon}_0}\log c_k'v_k'}dx\geq a_2>0,
$$
which implies that
$$
\int_{B_R}|c_k'v_k'|^{-\frac{\alpha}{\hat{\epsilon}_0}}dx<C,
\s \int_{B_R}|c_k'v_k'|^{\frac{\alpha}{\hat{\epsilon}_0}}dx<C.
$$
Recall that $c_k'v_k'$ satisfies
$$
\begin{array}{r}
-\Delta_{g_k(y_k+\rho(y_k)(x_k+r_kx))}c_k'v_k'+c(n)r_{k}^{2}R(g_k(y_k+\rho(y_k)(x_k+r_kx)))c_k'v_k'\\
=
c(n)R_k(y_k+\rho(y_k)(x_k+r_kx)){v_k'}^{\frac{4}{n-2}}c_k'v_k',
\end{array}
$$
and
$$
\int|R_k'|^p{v_k'}^\frac{2n}{n-2}dx<\Lambda.
$$
Choosing $\hat{\epsilon}_0$ to
be sufficiently small. By Lemma \ref{A1},
$c_k'v_k'$ is bounded in $W^{2,p}(B_R)$.
By volume identity, $v_k'\rightarrow 0$, then $R_k(g_k'){v_k'}^{\frac{4}{n-2}}$
converges to 0 in $L^q$ for any $q<p$. Thus $e^\phi$ is harmonic on $\R^n$ and
$$
\int_{B_1}|\nabla\phi|^2dx=\hat{\epsilon}_0^2.
$$
However, since $e^\phi$ is positive, $e^\phi$ is a constant which implies that
$$
\int_{B_1}|\nabla\phi|^2dx=0.
$$
We get a contradiction. \endproof

Now, we go back to the proof of Theorem \ref{main3}.
 We fix a  $B^{\R^n}_{r'}(x')\subset\subset\R^n\setminus\tilde{\CS}$ and
take $\tilde{c}_k$ such that

$$
\int_{B_{r'}(x')}\log \tilde{c}_k\tilde{v}_k=0.
$$
Since $\tilde{v}_k$ converges to $0$ weakly in $W^{2,p}(\R^n\setminus\tilde{\CS})$, $\tilde{c}_k\rightarrow+\infty$.
By Lemma \ref{brk}, $\|\nabla\log \tilde{v}_k\|_{L^2(\Omega)}$ is bounded for any
$\Omega\subset\subset\R^n\setminus\title{\CS}$.
By Poincar\'e inequality and Sobolev embedding,
we may assume $\log \tilde{c}_k\tilde{v}_k$ converges in $L^q(\Omega)$ for any $q<\frac{2n}{n-2}$,
and for a.e. $x\in\R^n\setminus\tilde{\CS}$.

Now, we fix an $\Omega\subset\subset \R^n\setminus\tilde{\CS}$ and
put
$$
\rho=\frac{1}{2}\min\{\inf_{x\in\Omega}\liminf_{k\rightarrow+\infty}r_k(x),\inf_{x\in\Omega}d(x,\tilde{\CS})\}
$$
which is positive. As \eqref{brk}, for any $y\in
\Omega$, we can prove $(\rho)^\frac{n-2}{2}\tilde{c}_{k}\tilde{v}_k(y+\rho x)$ converges to
a positive harmonic function weakly in $W^{2,p}(B_1)$.
Therefore, $\tilde{c}_k\tilde{v}_k$  converges to a positive harmonic function $\tilde{v}$
weakly in $W^{2,p}_{loc}(\R^n\setminus\tilde{\CS})$.
By Bocher's Theorem (see Theorem 3.9 in \cite{Sheldon-Paul-Wade}), $\tilde{v}(x)=v_0+a|x|^{2-n}$ near 0,
where $v_0$ is a smooth harmonic function defined on a neighborhood of $0$. By
\eqref{tk} and \eqref{concentrationonGamma},
$$
diam_{\tilde{v}^\frac{4}{n-2}g_{\R^n}}\Gamma_1
=\inf_{t\in(0,+\infty)}diam_{\tilde{v}^\frac{4}{n-2}g_{\R^n}}\Gamma_t.
$$
Hence $a\neq 0$. Then it is easy to check that $\tilde{v}^\frac{4}{n-2}g_{\R^n}$
is complete and noncompact near $0$. By Kelvin transformation, we can also check that
$\tilde{v}^\frac{4}{n-2}g_{\R^n}$ is complete and noncompact near $\infty$.
Also by Bocher's Theorem in \cite{Sheldon-Paul-Wade}, $\tilde{v}^\frac{4}{n-2}g_{\R^n}$ is complete
near each $y\in\tilde{\CS}$.\\

\noindent{\bf Case 2.}  $u_0> 0$.

For this case, if we replace $\hat c_k$ with $1$ and let $g_0=u_0^\frac{4}{n-2}g$ in the
proof for case 1, we complete the proof.\\

\noindent{\bf Case 3.} $\{u_k\}$ has at least 2 bubbles.

In this subsection, we assume $\{u_k\}$ has at least 2
bubbles. If $\{u_k\}$ has at least two concentration points, then
$G$ has at least 2 singular points. Thus, $G^\frac{4}{n-2}g$ is already
a complete manifold with $R=0$  which has at least 2 ends.
Therefore, we assume $\{u_k\}$ has only one concentration point $x_0$.
Let $(x_k^1,r_k^1)$ be a treetop. We have 2 subcases.

Subcase 1: there is a blowup sequence $(x_k^2,r_k^2)$ which is under
$(x_k^1,r_k^1)$. For this case, $v_k^2$ converges to nonzero function
which is in $W^{2,p}$ in a neighborhood of $0$, and $v^1$
can be considered as a bubble of $\{v_k^2\}$. Using the arguments in
the above subsection, we complete the proof.

Subcase 2: There is no blowup sequence $(x_k^2,r_k^2)$ under $(x_k^1,r_k^1)$.
Then we have $\frac{x_k^2-x_k^1}{r_k^1+r_k^2}\rightarrow+\infty$.
We set  $v_k^{2'}=|x_k^2-x_k^1|^{\frac{n-2}{2}}u_k(x_k^1+|x_k^2-x_k^1|x)$. Then both $v^1$
and $v^2$ can be considered as bubbles of $v_k^{2'}$ at 0. Since
$(x_k^1,|x_k^2-x_k^1|)$ is not a blowup sequence (otherwise, it
is a blowup sequence under $(x_k^1,r_k^1)$), $v_k^{2'}$ converges to 0
weakly in  $W^{2,p}(\R^n\setminus\CS')$, where $\CS'$ is the concentration
point set of $v_k^{2'}$. Then similar to Case 1, we can find $c_k$,
such that $c_kv_k^{2'}$ converges to a positive function $v^{2'}$
which is harmonic on $\R^n\setminus\CS'$. Let $\frac{x_k^2-x_k^1}
{|x_k^2-x_k^1|}\rightarrow x_0'$. Obviously, $0,x_0'\in \CS'$.
If $v^{2'}$ is smooth at $0$  or $x_0'$, we can use the arguments in
case 1 to complete the proof. If $v^{2'}$ is not smooth
at both $0$ and $x_0'$, then $(\R^n,(v^{2'})^\frac{4}{n-2}g_{\R^n})$
is a complete manifold with at least 2 ends.

\section{A Manifold sequence with $\int_{M}|Ric(g_k)|^p<\Lambda$}
Now, we let $g_k=u_k^\frac{4}{n-2}g$, and
$$
\vol(M,g_k)=1,\s \int_M|Ric(g_k)|^p\leq\Lambda.
$$
Assuming $g_k$ does not  converge weakly in $W^{2,p}(M,g)$.

If the Gromov-Hausdorff limit of $(M,g_k)$ is not a manifold $(S^n,g_\infty)
\in SC^p(\mathbb{S}^n)$, then one of the  followings must  happen:

1) $\{u_k\}$ has at least 2 bubbles at a point in $\CS$;
2) $u_0\neq 0$;
3) There exist at least 2 points in $\CS$, such that $\lambda_y>0$,
where $\lambda_y$ is as in Theorem \ref{main2}.

For all cases,
 we can find $c_k\rightarrow+\infty$,$t_{k}\rightarrow 0$, $y_k\in M$, such that
$(M,c_{k}^{\frac{4}{n-2}}t_{k}^{2}g_k(y_{k}+t_{k}x),y_k)$ converges to a $(M_\infty, g_\infty, y_\infty)$, and

i) $M_\infty$ is a smooth complete manifold with at least 2 ends, and each end is
collared topologically by $S^{n-1}\times\R$. $g_\infty$ is a smooth metric over $M_\infty$;

ii) $(M_\infty,g_\infty)$
is conformal to $(M,g)$ or $\mathbb{S}^{n}$ with finitely many points removed;

iii) $Ric(g_\infty)=0$.

\noindent{\it Proof of Corollary \ref{main4}:}
 By The Splitting Theorem, $(M_{\infty},g_{\infty})$ is isometric to $(N \times R^{k})$. If $k\geq 2$,  then given any compact subset $A$ of
$M_\infty$, $M_{\infty}\setminus A$ is connected ,which contradict the assumption $M_{\infty}$  has at least 2 ends.
Thus, $k=1$ and $M_\infty$
has exactly 2 ends.
Noting that  each end is
collared topologically by $S^{n-1}\times\R$, $M_\infty$
is diffeomorphic to $S^n$ with two points removed.
We have two cases: 1) $(M_\infty,g_\infty)$ is
conformal to $\R^n\setminus\{0\}$, 2) $M$ is diffeomorphic to $S^n$
and $(M_\infty,g_\infty)$ is conformal to $(M\setminus\{p_1,p_2\},g)$,
where $p_1$, $p_2\in M$.

When  $(M_\infty,g_\infty)$ is
conformal to $\R^n\setminus\{0\}$, we can set
$g_\infty=(ar^{2-n}+b)^\frac{4}{n-2}g_{\R^n}$, where $a$, $b>0$. Then we have
$$
\vol(B_1\setminus B_t(0),g_\infty)=\int_{B_1\setminus B_t}(ar^{2-n}+b)^\frac{2n}{n-2}dx\geq c_1 t^{-n},
$$
where $t$ is sufficiently small and $c_1>0$.
Choose a point $P_{0}$ on the $\partial B_{1}$. It is easy to check that
we can find $c_2>0$, such that for any $x\in\partial B_1$,
$$
d_{g_{\infty}}(x,sx)\leq c_{2}s^{-1},\s \forall s\in [t,1].
$$
Since, for any $y\in B_1\setminus B_t$,
$$
d_{g_\infty}(P_0,y)\leq diam_{g_\infty}(\partial B_1)+d_{g_\infty}(\frac{y}{|y|},y))\leq diam_{g_\infty}(\partial B_1)+c_{2}t^{-1},
$$
we can find $c_3>0$, such that $d_{g_\infty}(P_0,y)\leq c_3t^{-1}$, i.e.,
$B_{1}\setminus B_{t}(0)\subset B^{g_\infty}_{\frac{c_3}{t}}(P_{0})$.
Thus, we get
$$
\frac{\vol (B^{g_\infty}_{\frac{c_4}{t}}(P_{0}))}{(\frac{c_4}{t})^{n}}\geq \frac{c_{1}}{c_4^n}>0.
$$
On one hand, since $(M_\infty,g_\infty)$ is splitting to $S^{n-1}\times\R$, we must have
$$
\vol (B^{g_\infty}_{\frac{c_4}{t}}(P_{0}))=O(\frac{1}{t}),\s as\s t\rightarrow 0.
$$
So we get a contradiction.

For the case 2, we have $g_\infty=G^\frac{4}{n-2}g$, where $G$
is a Green function which satisfies
$$
-\Delta G+R_gG=\lambda_1p_1+\lambda_2p_2,
$$
where $\lambda_1$, $\lambda_2>0$. Then we have
$G=a_ir^{2-n}(1+o(1))$ in a normal  neighborhood  of $p_i$.
Then we can get a contradiction as in the case 1.

}


\end{document}